\documentclass[12pt]{article}
\title{Hyperplane arrangements and Lefschetz's hyperplane section theorem}
\author{Masahiko Yoshinaga \\
International Centre for Theoretical Physics\\
Strada Costiera 11, Trieste 34014, Italy\\
email: myoshina@ictp.it}
\date{December 11, 2006}
\usepackage{amssymb}

\newtheorem{Def}{Definition}[subsection]
\newtheorem{Prop}[Def]{Propositon}
\newtheorem{Thm}[Def]{Theorem}
\newtheorem{Lemma}[Def]{Lemma}
\newtheorem{Cor}[Def]{Corollary}

\newtheorem{Example}[Def]{Example}

\newcommand{\frakf}{\mathfrak{f}}
\newcommand{\frakg}{\mathfrak{g}}

\newcommand{\bbC}{\mathbb{C}}

\newcommand{\bbK}{\mathbb{K}}
\newcommand{\bbP}{\mathbb{P}}

\newcommand{\bbR}{\mathbb{R}}

\newcommand{\bbZ}{\mathbb{Z}}

\newcommand{\rmD}{\mathrm{D}}

\newcommand{\rmS}{\mathrm{S}}
\newcommand{\rmT}{\mathrm{T}}

\newcommand{\calA}{\mathcal{A}}
\newcommand{\calC}{\mathcal{C}}

\newcommand{\calF}{\mathcal{F}}
\newcommand{\calG}{\mathcal{G}}

\newcommand{\calL}{\mathcal{L}}

\newcommand{\calN}{\mathcal{N}}
\newcommand{\calO}{\mathcal{O}}
\newcommand{\calP}{\mathcal{P}}
\newcommand{\calT}{\mathcal{T}}
\newcommand{\calU}{\mathcal{U}}
\newcommand{\calX}{\mathcal{X}}

\newcommand{\bch}{\mathsf{bch}}

\newcommand{\ch}{\mathsf{ch}}
\newcommand{\cpl}{\mathsf{M}}

\newcommand{\sfch}{\mathsf{ch}}

\newcommand{\rmCrit}{\mathrm{Crit}}
\newcommand{\rmCyl}{\mathrm{Cyl}}
\newcommand{\rmGal}{\mathrm{Gal}}

\newcommand{\rmint}{\mathrm{int}}
\newcommand{\M}{\mathsf{M}}
\newcommand{\sfM}{\mathsf{M}}

\newcommand{\rmPr}{\mathrm{Pr}}
\newcommand{\Poin}{\mathrm{Poin}}
\newcommand{\rmVect}{\mathrm{Vect}}

\newcommand{\ii}{\sqrt{-1}}

\makeatletter
\def\codim{\mathop{\operator@font codim}\nolimits}
\def\grad{\mathop{\operator@font grad}\nolimits}
\def\Hom{\mathop{\operator@font Hom}\nolimits}
\def\Im{\mathop{\operator@font Im}\nolimits}
\def\Ker{\mathop{\operator@font Ker}\nolimits}
\def\Loc{\mathop{\operator@font Loc}\nolimits}
\def\rank{\mathop{\operator@font rank}\nolimits}
\def\Rep{\mathop{\operator@font Rep}\nolimits}
\def\tr{\mathop{\operator@font tr}\nolimits}
\def\sign{\mathop{\operator@font sign}\nolimits}
\makeatother

\newcommand{\owari}{\hfill$\square$}

\begin{document}
\maketitle

\begin{abstract}
The Lefschetz hyperplane section theorem 
asserts that a complex affine variety is homotopy equivalent to a 
space obtained from its generic hyperplane section by attaching 
some cells. 
The purpose of this paper is to give an explicit description of 
attaching maps of these cells for the complement of a complex hyperplane 
arrangement defined over real numbers. 
The cells and attaching maps are described in combinatorial 
terms of chambers. 
We also discuss the 
cellular chain complex with coefficients in a local system and 
a presentation for the fundamental group associated to 
the minimal CW-decomposition for the complement. \\
{\bf MSC-class}: 32S22 (Primary) 14F35, 32S50 (Secondary)\\
{\bf Keywords}: Arrangements, Minimal CW-complex, Homotopy type. 
\end{abstract}

\section{Introduction}
\label{sec:intro}

The Lefschetz hyperplane section theorem is a result 
concerning a topological relationship between 
an algebraic variety and its generic hyperplane section. 
%A tremendous amount of effort has gone into generalizing 
%this theorem. 
The following is a version of the Lefschetz theorem 
for affine varieties. 
Let $g\in\bbC[x_1, \ldots, x_\ell]$ 
be a polynomial and $\cpl(g):=\{x\in\bbC^\ell|\ g(x)\neq 0\}$ be the 
hypersurface complement defined by $g$.

\begin{Thm}
\label{thm:zf}
{\normalfont\bf (Affine Lefschetz Theorem \cite{ham-aff, le-ham})}
Let $F$ be a generic affine hyperplane in $\bbC^\ell$. Then 
the space $\cpl(g)$ has the homotopy type of 
a space obtained from $\cpl(g)\cap F$ by attaching 
a certain number of $\ell$-dimensional cells. 
\end{Thm}
The important part of the above Lefschetz theorem for affine varieties 
is that the cells attached to $\cpl(g)\cap F$ 
all have equal dimension $\ell$. 
This makes the situation relatively simple. 
An immediate corollary, obtained by induction on the dimension $\ell$, 
is that $\cpl(g)$ is homotopy equivalent to an 
$\ell$-dimensional CW-complex whose 
$(\ell-1)$-skeleton is homotopy equivalent to 
$\cpl(g)\cap F$, and we also conclude that 
the number of $\ell$-cells is equal to 
$\dim H_\ell(\cpl(g), \cpl(g)\cap F)$. 
The number of $\ell$-cells is obviously greater than or equal to 
the Betti number $b_\ell(\cpl(g))$. More precisely, we have 
the following exact sequence: 
$$
0\rightarrow 
H_\ell(\cpl(g))\rightarrow 
H_\ell(\cpl(g), \cpl(g)\cap F)\rightarrow 
H_{\ell-1}(\cpl(g)\cap F)
\stackrel{i_{\ell-1}}{\longrightarrow}
H_{\ell-1}(\cpl(g)).
$$

Another corollary is 
\begin{Cor}
Let $i_p:H_p(\cpl(g)\cap F, \bbC)\rightarrow H_p(\cpl(g), \bbC)$ denote the 
homomorphism induced from the natural inclusion 
$i:\cpl(g)\cap F\hookrightarrow \cpl(g)$, then 
$$
i_p\mbox{\normalfont{ is }} 
\left\{
\begin{array}{ll}
\mbox{\normalfont{ isomorphic }} &\mbox{\normalfont{ for }}
p=0, 1, \ldots, \ell-2\\
\mbox{\normalfont{ surjective }} &\mbox{\normalfont{ for }}
p=\ell-1. 
\end{array}
\right. 
$$ 
\end{Cor}

As noted by A. Dimca, S. Papadima and R. Randell  
(\cite{dim-pap}, \cite{ran-mor}), suppose $i_{\ell-1}$ is isomorphic, 
then the number of $\ell$-dimensional cells attached would be equal to 
the Betti number $b_\ell(\cpl(g))$. 
While in case of a hyperplane arrangement, 
that is, when $g$ is a product of linear equations, 
$i_{\ell-1}$ is indeed isomorphic 
(see Prop. \ref{prop:trunc}), and hence the number 
of $\ell$-cells is exactly equal to $b_\ell(\cpl(g))$. 
%Thus in contrast with the situation of general varieties, 
%complements of hyperplane arrangements behave very 
%nicely from the Lefschetz theorem point of view. 

Repeating the same procedure inductively, we finally obtain a 
minimal CW decomposition. 

\begin{Thm}
{\normalfont\bf (\cite{dim-pap} \cite{ran-mor})}
Let $\calA$ be an affine arrangement in $\bbC^\ell$. Then 
the complement $\sfM(\calA)$ is homotopy equivalent to 
a minimal {\rm CW}-complex, i.e. a {\rm CW}-complex 
whose number of $k$-cells is equal to 
$b_k(\sfM(\calA))$ for each $k$. 
\end{Thm}

Let $\calL$ be a rank one local system on $\sfM(\calA)$. 
Then the minimal CW-decomposition yields 
a cellular chain complex 
$(\calC_\bullet(\sfM(\calA), \calL), \partial)$ satisfying 
$\dim\calC_k(\sfM(\calA), \calL)=b_k(\sfM(\calA))$ and 
$H_k(\sfM(\calA), \calL)\cong H_k(\calC_\bullet(\sfM(\calA), \calL), \partial)$. 
We call this the twisted minimal chain complex. 
This kind of minimal complexes were first constructed by 
D. Cohen by using stratified Morse theory \cite{coh-int}. 
Properties of twisted minimal chain complexes 
have been studied in many papers including 
\cite{coh-mor, coh-orl, dp-equiv, ps-h}. 
To describe boundary maps 
$\partial:\calC_\bullet\rightarrow\calC_{\bullet-1}$, 
some information about the attaching maps of minimal CW-complexes 
are required. 
Attaching maps for 
minimal CW-decompositions for $\ell=2$ were studied by 
M. Falk \cite{fal-hom} based on \cite{ran-fun, sal-top} 
(see also \cite{lib}). 

However, little is known about both the attaching maps 
and the boundary maps 
$\partial:\calC_\bullet\rightarrow\calC_{\bullet-1}$ for 
higher dimensional cases. 
The purpose of this paper is to describe how 
$\ell$-cells are attached to a generic hyperplane 
section $\sfM(\calA)\cap F$ for the complement $\sfM(\calA)$ of 
a real hyperplane arrangement $\calA$ (\S\ref{subsec:const}). Here, 
``real hyperplane arrangement'' 
means that the defining polynomial 
$g\in\bbR[x_1, \ldots, x_\ell]$ is a product of linear 
equations with real coefficients. 
Although we have not yet obtained a complete understanding 
of the minimal CW-decomposition for hyperplane complements, 
we obtain a description of the twisted minimal chain complex 
$(\calC_\bullet(\sfM(\calA), \calL), \partial)$ (in \S\ref{sec:twist}). 
Our formula of the twisted boundary map contains the integers 
``$\deg(C, C')$''. It is defined by using the topological 
relationship between two chambers, and its computation will 
be quite difficult in general. 
But it is computable in a certain cases. 
Our presentation has some 
applications on the structure of local system homologies, which will be 
discussed in a subsequent paper \cite{yos-res}. 

%But unfortunately, its computation 
%will be quite difficult. 
%Thus I do not know whether if our description 
%is applicable to computations. 

The advantage of focusing our attention on real arrangements 
is that we can use 
structures of chambers, namely, the connected 
components of $\sfM(\calA)\cap \bbR^\ell$. 
The study of relationships between topology of $\sfM(\calA)$ 
and combinatorics of chambers is a classical topic in 
the theory of hyperplane arrangements. 
We summarize some classical results related to 
chamber-counting problems in \S\ref{sec:comb}. 
The number of chambers are related to Betti numbers of 
$\sfM(\calA)$. 
Later, in \S\ref{subsec:stable}, we give a more geometric interpretation to 
these numerical relations between chambers and Betti numbers: 
chambers can be thought of as stable manifolds 
for a certain Morse function. 
This interpretation will play a crucial role in this paper. 
By a well-known duality between stable and unstable manifolds, 
the set of chambers are indexing unstable cells which appear in the 
minimal CW-decomposition. Thus the basis of the associated cellular 
chain complex is also indexed by chambers.

In \S\ref{sec:sd} we review the Salvetti complex and 
the Deligne groupoid. They relate combinatorial structures of 
chambers to topological structures of the complexified 
complements. In particular, for the purposes of this paper, 
we have to describe local systems in terms of chambers. 
The Deligne groupoid offers 
an appropriate language to deal with 
local systems in a combinatorial context. 
A local system can be interpreted 
as a representation of the Deligne groupoid.

In \S\ref{subsec:lef} we give a proof of the 
Theorem \ref{thm:zf} for 
hyperplane complements. It is proved by applying Morse theory 
to a Morse function of the form $|f/g|^2$, where $f$ is a defining 
equation of the generic hyperplane $F$. 
Although the proof does not 
involve anything new, Morse theoretic consideration 
in the proof will be needed later. 
In particular, Morse theory tells us that, 
under Morse-Smale condition on the gradient vector field, 
the unstable manifolds can be viewed as 
the $\ell$-cells attached to the generic section, 
and we have a homotopy 
equivalence 
$$
\sfM(\calA)\approx
(\sfM(\calA)\cap F)\cup\bigcup_{p\in\rmCrit(\varphi)}W_p^u, 
$$
where $W_p^u$ is the unstable manifold corresponding to 
a critical point $p\in\rmCrit(\varphi)$ of the Morse function $\varphi$. 
From the Morse-Smale condition, unstable and stable manifolds 
define ``set-theoretic'' dual 
bases of $H_\ell(\sfM(\calA))$ and $H_\ell^{lf}(\sfM(\calA))$, respectively, 
that is, 
\begin{equation}
\label{eq:set}
W_p^u\cap W_q^s=
\left\{
\begin{array}{cl}
W_p^u\pitchfork W_q^s=\{p\}&\mbox{ if }p=q, \\
\emptyset & \mbox{ if } p\neq q. 
\end{array}
\right.
\end{equation}
The main result in \S\ref{sec:morse} is that the set-theoretical 
duality between stable and unstable manifolds characterizes the homotopy 
type of unstable manifolds. 

As noted above, in the case of real arrangements, a stable manifold 
is known to be equal to a chamber. The goal of \S\ref{sec:constcell} 
is to construct cells attached to the hyperplane section 
$\sfM(\calA)\cap F$ which satisfy the set-theoretical duality condition 
(\ref{eq:set}) with respect to the chambers. 
Thanks to the result in the previous section, 
the cells constructed in this section are homotopy equivalent to 
unstable manifolds. The special case when $\ell=2$ offers 
a new presentation for the fundamental group $\pi_1(\sfM(\calA))$, 
which is given in the appendix \S\ref{sec:fundam}. 

In \S\ref{sec:twist}, using the construction of the 
cells in the previous section, 
we determine the boundary map of twisted cellular complex 
of the minimal CW-decomposition. The essential ingredient is 
calculating twisted intersection numbers of the boundary of 
a cell and chambers. In \S\ref{subsec:deg}, we introduce the 
concept of the degree map which associates to a pair of chambers 
$(C, C')$ an integer $\deg(C, C')$. The degree map is required 
for both the boundary maps of twisted cellular chain complexes 
and the presentations for fundamental groups.

\section{Combinatorics of arrangements}
\label{sec:comb}

In this section we establish some relationships 
among generic subspaces, 
numbers of chambers and Betti numbers for 
complements of hyperplane arrangements. 

\subsection{Basic constructions}
\label{sec:basic}

Let $V$ be an $\ell$-dimensional vector space. 
A finite set of affine hyperplanes 
$\calA=\{H_1, \ldots, H_n\}$ is called a {\it hyperplane arrangement}. 
Let $L(\calA)$ be the set of nonempty intersections of 
elements of $\calA$. Define a partial order on $L(\calA)$ 
by $X\leq Y\Longleftrightarrow Y\subseteq X$ for $X, Y\in L(\calA)$. 
Note that this is reverse inclusion. 

Define a {\it rank function} on $L(\calA)$ by $r(X)=\codim X$. 
Write $L^p(\calA)=\{X\in L(\calA)|\ r(X)=p\}$. We call $\calA$ 
{\it essential} if $L^\ell(\calA)\neq\emptyset$. 

Let $\mu:L(\calA)\rightarrow \bbZ$ be the {\it M\"obius function} 
of $L(\calA)$ defined by 
$$
\mu(X)=
\left\{
\begin{array}{ll}
1 &\mbox{ for }X=V\\
-\sum_{Y<X}\mu(Y), &\mbox{ for } X>V. 
\end{array}
\right.
$$
The {\it Poincar\'e polynomial} of $\calA$ is 
$\pi(\calA, t)=\sum\nolimits_{X\in L(\calA)}\mu(X)(-t)^{r(X)}$ and 
we also define numbers $b_i(\calA)$ by the formula 
$$
\pi(\calA, t)
=\sum\nolimits_{i=0}^\ell b_i(\calA)t^i. 
$$
We also define the {\it $\beta$-invariant} $\beta(\calA)$ by 
$$
\beta(\calA)=|\pi(\calA, -1)|, 
$$
if $\calA$ is an essential arrangement, the sign can be precisely 
enumerated as $\beta(\calA)=(-1)^\ell\pi(\calA, -1)$. 

Given a hyperplane $H\in\calA$, we define other arrangements: 
the {\it deletion} of $\calA$ with respect to $H$ is 
$\calA'=\calA\setminus\{H\}$ and the {\it restriction} is 
$\calA''=\{H\cap K\ |\ K\in\calA'\}$. Note that the restriction 
$\calA''$ is an arrangement in $H$. 
The Poincar\'e polynomials satisfy the following recursion: 
\begin{equation}
\label{eq:rec}
\pi(\calA, t)=\pi(\calA', t)+t\cdot\pi(\calA'', t). 
\end{equation}

\subsection{Classical results}
\label{sec:classic}

Let $\calA$ be an arrangement in a real vector space 
$V_\bbR$. 
Then the complement $V_\bbR\setminus \bigcup_{i=1}^n H_i$ is a union 
of open, connected components called {\it chambers}. 
Let us denote the set of all chambers by $\ch(\calA)$, and 
the set of relatively compact (or bounded) chambers by $\bch(\calA)$. 
If $\calA$ is an arrangement in a complex vector space 
$V_\bbC$, the complement is a connected affine algebraic 
variety and denoted by $\M(\calA)$.  

The Poincar\'e polynomial defined above and the geometric structure of 
the complement are related by the following theorem. 

\begin{Thm}{\normalfont \cite{orl-sol, zas-face}}
\label{thm:classic}
\begin{itemize}
\item[{\normalfont (i) }] 
Let $\calA$ be an essential real $\ell$-arrangement. 
The number $|\ch(\calA)|$ of chambers and the number 
$|\bch(\calA)|$ of bounded chambers $|\bch(\calA)|$ are given by 
\begin{eqnarray*}
|\ch(\calA)|&=&\pi(\calA, 1)\\
|\bch(\calA)|&=&(-1)^\ell \pi(\calA, -1)=\beta(\calA). 
\end{eqnarray*}
\item[{\normalfont (ii) }] 
Let $\calA$ be a complex arrangement. 
Then $b_i(\calA)$ is equal to the topological Betti number 
$b_i(\sfM(\calA))$. In other words, 
the topological Poincar\'e polynomial 
$\Poin(\sfM(\calA), t)=\sum_i b_i(\sfM(\calA)) t^i$
is given 
by 
$$
\Poin(\sfM(\calA), t)=\pi(\calA, t). 
$$
In particular, the absolute value of the 
topological Euler characteristic $|\chi(\sfM(\calA))|$ 
of the complement is equal to $\beta(\calA)$. 
\end{itemize}
\end{Thm}

\subsection{Generic flags}
\label{sec:flag}

Let $\calA$ be an $\ell$-arrangement. 
A $q$-dimensional affine subspace $\calF^q\subset V$ is 
called {\it generic} or {\it transversal} to $\calA$ if 
$\dim \calF^q\cap X=q-r(X)$ for $X\in L(\calA)$. 
A {\it generic flag} 
$\calF$ is defined to be a complete flag (of affine subspaces) in $V$, 
$$
\calF:\ \emptyset=\calF^{-1}\subset\calF^0\subset\calF^1\subset\cdots 
\subset\calF^\ell=V,  
$$
where each $\calF^q$ is a generic $q$-dimensional affine subspace. 

For a generic subspace $\calF^q$ we have an arrangement in $\calF^q$ 
$$
\calA\cap\calF^q:=\{H\cap\calF^q|\ H\in\calA\}. 
$$
The genericity provides an isomorphism of posets 
\begin{equation}
\label{eq:poset}
L(\calA\cap\calF^q)\cong L^{\leq q}(\calA):=\bigcup_{i\leq q}L^i(\calA). 
\end{equation}
In \cite{orl-sol} Orlik and Solomon gave a presentation 
of the cohomology ring $H^*(\sfM(\calA), \bbZ)$ in terms 
of the poset $L(\calA)$ for a complex arrangement $\calA$. 
The next proposition follows from (\ref{eq:poset}). 
\begin{Prop}
\label{prop:trunc}
Let $\calA$ be a complex arrangement and $\calF^q$ a 
$q$-dimensional generic subspace. Then the natural inclusion 
$i:\sfM(\calA)\cap\calF^q\hookrightarrow \sfM(\calA)$ 
induces isomorphisms 
$$
i_k: H_k(\sfM(\calA)\cap\calF^q, \bbZ)\stackrel{\cong}{\longrightarrow}
H_k(\sfM(\calA), \bbZ), 
$$
for $k=0, 1, \ldots, q$. 
\end{Prop}
In particular, the Poincar\'e polynomial of $\calA\cap\calF^q$ 
is given by 
\begin{equation}
\label{eq:trun}
\pi(\calA\cap\calF^q, t)=\pi(\calA, t)^{\leq q}, 
\end{equation}
where $(\sum_{i\geq 0}a_i t^i)^{\leq q}=\sum_{i=0}^q a_i t^i$ is the 
truncated polynomial. These formulas and 
Theorem \ref{thm:classic} prove the following result. 

\begin{Prop}
\label{prop:flag}
Let $\calA$ be a real $\ell$-arrangement and $\calF$ a generic flag. 
Define 
$$
\sfch_q^\calF(\calA)=\{C\in\sfch(\calA)|\ C\cap\calF^q\neq\emptyset 
\mbox{ {\normalfont and} }C\cap\calF^{q-1}=\emptyset \}, 
$$
for each $q=0, 1, \ldots, \ell$. 
Then 
\begin{itemize}
\item[{\normalfont (i)}] $|\sfch_q^\calF(\calA)|=b_q(\sfM(\calA)) $. 
\item[{\normalfont (ii)}] If $\calA$ is essential, then 
$b_\ell(\sfM(\calA))=\beta(\calA\cup \{\calF_{\ell-1}\})$, 
\end{itemize}
where $\sfM(\calA)$ is the complement of the complexified arrangement 
of $\calA$ and $\calA\cup\{\calF_{\ell-1}\}$ is the arrangement 
obtained by adding $\calF^{\ell-1}$ to $\calA$. 
%added with a hyperplane $\calF^{\ell-1}$. 
\end{Prop}

\noindent
{\bf Proof.} For any chamber $C\in\sfch(\calA)$, the intersection 
$C\cap \calF^q$ is either an empty set or a chamber in $\calA\cap\calF^q$. 
Hence we have a bijection 
$$
\begin{array}{ccc}
\bigcup_{i\leq q}\sfch_i^\calF(\calA)&\longrightarrow &\sfch(\calA\cap\calF^q)\\
&&\\
C&\longmapsto& C\cap\calF^q. 
\end{array}
$$
Counting the number of chambers by using Theorem \ref{thm:classic} (i) 
and (\ref{eq:trun}), we obtain 
\begin{eqnarray*}
\sum_{i\leq q}|\sfch_q^\calF(\calA)|&=&|\sfch(\calA\cap\calF^q)|\\
&=&\pi(\calA\cap\calF^q, t)|_{t=1}\\
&=&\sum_{i\leq q}b_i(\calA). 
\end{eqnarray*}
Thus we have (i). 

The recursion formula (\ref{eq:rec}) allows us to calculate 
the Poincar\'e polynomial of $\calA\cup\calF^{\ell-1}$: 
\begin{eqnarray*}
\pi(\calA\cup\{\calF^{\ell-1}\}, t)&=&\pi(\calA, t)+t\cdot\pi(\calA\cap\calF^{\ell-1}, t)\\
&=&\pi(\calA, t)+t\cdot\pi(\calA, t)^{\leq\ell-1}\\
&=&\pi(\calA, t)+t\cdot\left(\pi(\calA, t)-b_\ell(\calA)t^\ell\right). 
\end{eqnarray*}
By putting $t=-1$ we obtain (ii). \owari

\begin{Example}\normalfont
{\normalfont Figure \ref{fig:flag}} shows an example of 
arrangement $\calA$ of three lines in $\bbR^2$ with a generic flag 
$\calF: \calF^0\subset\calF^1$. Note that $\pi(\calA, t)=1+3t+3t^2$. 

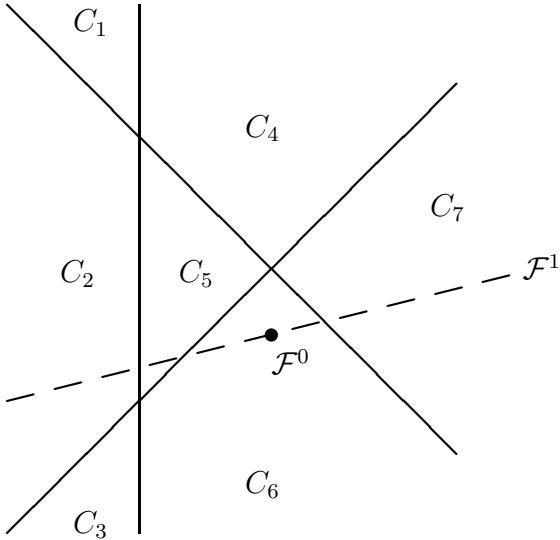
\begin{figure}[htbp]
\begin{picture}(100,190)(0,0)
\thicklines

\put(0,0){\line(1,1){170}}
\put(0,200){\line(1,-1){170}}
\put(50,0){\line(0,1){200}}
\multiput(0,50)(18,4.5){11}{\line(4,1){10}}
\put(100,75){\circle*{5}}
\put(25,190){$C_1$}
\put(20,95){$C_2$}
\put(25,0){$C_3$}
\put(90,150){$C_4$}
\put(65,95){$C_5$}
\put(90,15){$C_6$}
\put(160,120){$C_7$}

\put(100,60){$\calF^0$}
\put(195,96){$\calF^1$}

\put(260,140){A generic flag}
\put(280,120){$\calF: \{{\mathrm{pt}}\}=\calF^0\subset\calF^1\subset V$}

\put(260,100){$\sfch(\calA)=\{C_1, \ldots, C_7\}$}
\put(260,80){$\sfch_0^\calF(\calA)=\{C_6\}$}
\put(260,60){$\sfch_1^\calF(\calA)=\{C_2, C_5, C_7\}$}
\put(260,40){$\sfch_2^\calF(\calA)=\{C_1, C_3, C_4\}$}

\end{picture}
     \caption{A $2$-arrangement and a generic flag}\label{fig:flag}
\end{figure}

\end{Example}

Let $\calA$ be a real arrangement with a generic flag $\calF$. 
Consider the $\ell$-th homology, cohomology and homology 
with locally finite chains for the complement. 
Both $H_\ell^{lf}(\sfM(\calA), \bbC)$ and 
$H^\ell(\sfM(\calA), \bbC)$ are dual to 
$H_\ell(\sfM(\calA), \bbC)$. So there exists a canonical isomorphism 
\begin{equation}
\label{eq:isom}
H_\ell^{lf}(\sfM(\calA), \bbC)
\stackrel{\cong}{\longrightarrow} H^\ell(\sfM(\calA), \bbC). 
\end{equation}
Let $C$ be a chamber. Using the inclusion 
$V_\bbR\hookrightarrow V_\bbC=V_\bbR\oplus \ii V_\bbR$, 
$C$ can be considered as a locally finite $\ell$-dimensional cycle in 
$\sfM(\calA)$ and 
determines an element 
$[C]\in H_\ell^{lf}(\sfM(\calA))$. 

Recall that $C\in\sfch_\ell^\calF(\calA)$ is a chamber satisfying 
$C\cap\calF^{\ell-1}=\emptyset$, and that the number of such chambers is 
equal to the $\ell$-th Betti number 
$b_\ell(\calA)=\dim H_\ell^{lf}(\sfM(\calA))$. 
Later we will prove that 
$\left\{[C]\left|\ C\in\sfch_\ell^\calF(\calA)\right. \right\}$ 
forms a basis of $H_\ell^{lf}(\sfM(\calA))$ (Cor. \ref{cor:basis}). 

%\begin{Problem}
%Suppose $C\in\sfch_\ell^\calF(\calA)$. 
%Describe the image of $[C]$ in $H^\ell(\sfM(\calA), \bbC)$ 
%under the isomorphism {\normalfont (\ref{eq:isom})}. 
%\end{Problem}

\section{The Salvetti complex and the Deligne groupoid}

\label{sec:sd}

In \cite{sal-top} Salvetti has given a finite regular 
CW-complex which carries the homotopy type of the complement 
$\sfM(\calA)$ in the case where $\calA$ is a complexified real 
arrangement. In this section we review some results on the 
complexified complement $\sfM(\calA)$ of 
a real arrangement $\calA$. 

\subsection{Complexified real arrangements}
\label{sec:cpxf}

Let $\calA_\bbR$ be an arrangement in a real vector space 
$V_\bbR$. By definition each hyperplane $H\in\calA_\bbR$ 
is defined by a real equation $\alpha_H=0$ of degree one. 
The complexification $\calA_\bbC$ is a set of hyperplanes 
in $V_\bbC=V_\bbR\otimes\bbC$ defined by real equations 
$\alpha_H=0$ for $H\in\calA_\bbR$. 

Since $V_\bbC\cong V_\bbR\oplus \ii V_\bbR$, $V_\bbC$ can be 
identified with the total space of the tangent bundle 
$\rmT V_\bbR$. More precisely we identify as follows: 
\begin{equation}
\label{eq:tg}
\begin{array}{cccc}
\rmT V_\bbR&\stackrel{\cong}\longrightarrow&V_\bbC\\
&&&\\
(x, v)&\longmapsto&(x, v)_\bbC&=x+\ii v, 
\end{array}
\end{equation}
where 
$\rmT V_\bbR=\{(x, v)|\ x\in V_\bbR, v\in \rmT_xV_\bbR\cong V_\bbR\}$. 
This identification (\ref{eq:tg}) enables us to express a point 
in $V_\bbC$ as a tangent vector on $V_\bbR$, and a path 
in $V_\bbC$ can be expressed 
as a continuous family of tangent vectors 
along a path in $V_\bbR$, for simplicity we say a vector field along a 
path in $V_\bbR$. 

\begin{Example}\normalfont
The left side of 
{\normalfont Figure \ref{fig:vf}} expresses a vector field along the 
segment $[-1, 1]$ in $V_\bbR\cong \bbR$. The right side expresses the 
corresponding path in $V_\bbC\cong \bbC$. 

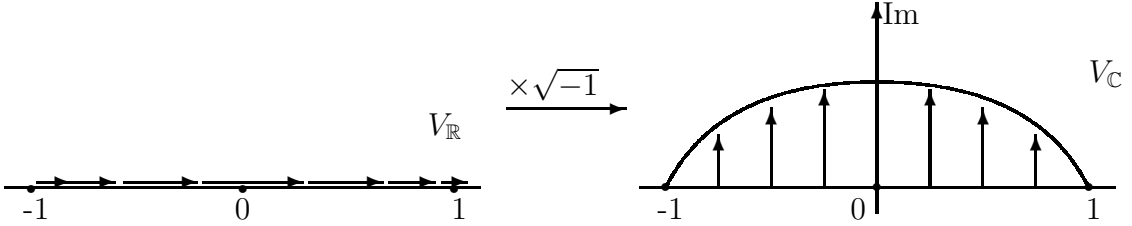
\begin{figure}[htbp]
\begin{picture}(100,80)(0,0)
\thicklines

\put(-10,0){\line(1,0){180}}
\put(230,0){\line(1,0){180}}

\put(180,30){\vector(1,0){45}}
\put(180,35){$\times\ii$}

\put(-3,-12){-1}
\put(77,-12){0}
\put(159,-12){1}
\multiput(0,-0.5)(80,0){3}{\circle*{3}}

\put(155,2){\vector(1,0){10}}
\put(135,2){\vector(1,0){18}}
\put(105,2){\vector(1,0){28}}
\put(65,2){\vector(1,0){38}}
\put(35,2){\vector(1,0){28}}
\put(14,2){\vector(1,0){18}}
\put(2,2){\vector(1,0){13}}

\put(237,-12){-1}
\put(310,-12){0}
\put(399,-12){1}
\multiput(240,0)(80,0){3}{\circle*{3}}

\qbezier(240,0)(260,40)(320,40)
\qbezier(400,0)(380,40)(320,40)

\multiput(260,0)(120,0){2}{\vector(0,1){20}}
\multiput(280,0)(80,0){2}{\vector(0,1){30}}
\multiput(300,0)(40,0){2}{\vector(0,1){37}}

\put(320,-10){\vector(0,1){80}}
\put(322,62){$\Im$}

\put(150,20){$V_\bbR$}
\put(400,40){$V_\bbC$}

\end{picture}
     \caption{Vector field along the segment $[-1,1]$ and corresponding path}
\label{fig:vf}
\end{figure}

\end{Example}

Let $x\in V_\bbR$. Then $\alpha_H(x)$ can be expressed as 
$\alpha_H(x)=a\cdot x+b$, where $a\in V_\bbR^*$  and 
$b\in\bbR$. Hence 
$$
\alpha_H(x+\ii v)=\alpha_H(x)+ \ii a\cdot v, 
$$
for $x+\ii v\in V_\bbC$. We have 
$$
\alpha_H(x+\ii v)=0\Longleftrightarrow \alpha_H(x)=0 \mbox{ and } a\cdot v=0. 
$$
This proves the following. 

\begin{Lemma}
\label{lem:compl}
Let $\calA$ be a real arrangement. For $x\in V_\bbR$ we define 
$\calA_x$ to be the set $\{H\in\calA | H\ni x\}$ of all hyperplanes 
containing $x$. Then the complexified complement is 
$$
\sfM(\calA)\cong\{(x, v)_\bbC| x\in V_\bbR,\ v\in\rmT_xV_\bbR\setminus \calA_x\}. 
$$
\end{Lemma}

\subsection{The Salvetti complex}

We recall some notions about the Salvetti complex, for details see 
\cite{om}. 

\begin{Def}\normalfont
Let $X\in L(\calA)$ be an intersection of a real arrangement $\calA$. 
A connected 
component $X^\circ$ of $X\setminus \bigcup_{H\nsupseteq X}H$ 
is called a {\it face} of $\calA$. The set of all faces is 
denoted by $\calL$. Define a partial order by 
$$
X\leq Y\Longleftrightarrow X\subset \bar{Y}, \mbox{ for }
X, Y\in\calL,  
$$
where $\bar{Y}$ is the closure of $Y$ in $V_\bbR$. 
The ordered set $(\calL, \leq)$ is called the {\it face poset} of $\calA$. 
\end{Def}

In this notation 
$\sfch(\calA)$ is the set of maximal elements in $(\calL, \leq)$. 

Given a face $X\in\calL$ and a chamber $C\in\sfch$, the 
chamber $X\circ C$ satisfying the following conditions is 
uniquely determined (see \cite{bhr} for more on $X\circ C$). 
%See some combinatorial application to bidigare hanlon someone 
%they are interesting
\begin{itemize}
\item[(1)]$X\leq X\circ C$, and 
\item[(2)]If $X$ is contained in a hyperplane 
$H\in\calA$, then $C$ and $X\circ C$ are 
on the same side with respect to $H$. 
\end{itemize}

\begin{Def}\normalfont
The poset $(\calP(\calA), \preceq)$ is defined as follows: 
\begin{eqnarray*}
&&\calP(\calA)=\{(X, C)\in\calL\times\sfch(\calA)\ |\ X\leq C\}\\
&&(X_1, C_1)\preceq(X_2, C_2)\Longleftrightarrow X_1\geq X_2 
\mbox{ and }X_1\circ C_2=C_1. 
\end{eqnarray*}
\end{Def}

\begin{Thm}
\label{thm:salvetti}
There exists a regular {\rm CW}-complex $X$, called the Salvetti complex, 
such that the face poset $\calF(X)$ of the complex $X$ is isomorphic to 
$\calP(\calA)$, and $X$ is homotopy equivalent to $\sfM(\calA)$. 
\end{Thm}

\begin{Example}\normalfont
We show some examples of low dimensional cells. 
\begin{itemize}
\item[{\normalfont ($0$-cell)}]
In $\calP(\calA)$, the $0$-cells of $X$ are corresponding to 
the $(C, C)\in\calP(\calA)$, $C\in\sfch(\calA)$. 

\item[{\normalfont ($1$-cell)}]
Two chambers $C$ and $C'$ are adjacent if $\bar{C}\cap\bar{C'}$ is 
contained in a hyperplane and has nonempty interior in the hyperplane. 
The relative interior of $\bar{C}\cap\bar{C'}$ is called the wall 
separating $C$ and $C'$. 
Let $C$ and $C'$ be adjacent chambers separated by a 
wall $X$. Then we have two $1$-cells, $(X, C)$ and $(X, C')$, 
which connect $(C, C)$ and $(C', C')$. 
{\normalfont (Figure \ref{fig:1-cell})}

\begin{figure}[htbp]
\begin{picture}(100,100)(0,0)
\thicklines

\put(50,0){\line(0,1){100}}
\put(30,0){\line(2,1){140}}
\put(30,100){\line(2,-1){140}}
\put(50,70){\circle*{5}}
\put(50,70){\line(-2,1){15}}
\put(20,80){$X$}
\put(65,60){$C$}
\put(20,25){$C'$}
\put(18,50){\circle*{3}}
\put(20,50){\vector(1,0){6}}
\put(27,50){\vector(1,0){10}}
\put(40,50){\vector(1,0){20}}
\put(62,50){\vector(1,0){10}}
\put(73,50){\vector(1,0){6}}
\put(80,50){\circle*{3}}

\put(250,0){\line(0,1){100}}
\put(230,0){\line(2,1){140}}
\put(230,100){\line(2,-1){140}}
\put(250,70){\circle*{5}}
\put(250,70){\line(-2,1){15}}
\put(220,80){$X$}
\put(265,60){$C$}
\put(220,25){$C'$}
\put(218,50){\circle*{3}}
\put(225,50){\vector(-1,0){6}}
\put(237,50){\vector(-1,0){10}}
\put(260,50){\vector(-1,0){20}}
\put(272,50){\vector(-1,0){10}}
\put(278,50){\vector(-1,0){6}}
\put(280,50){\circle*{3}}

\end{picture}
     \caption{$1$-cells corresponding to $(X, C)$ and $(X, C')$}\label{fig:1-cell}
\end{figure}
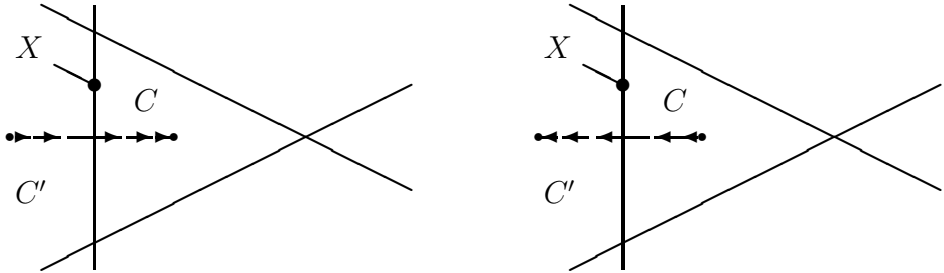

\item[{\normalfont ($2$-cell)}]
Let $X\in\calL$ be a face of codimension two with a chamber $C_1\geq X$. 
We have a $2$-cell $(X, C_1)$. 
{\normalfont (Figure \ref{fig:2-cell})}

\begin{figure}[htbp]
\begin{picture}(100,195)(0,0)
\thicklines

\put(80,100){\line(1,0){240}}
\put(150,0){\line(1,2){100}}
\put(150,200){\line(1,-2){100}}
\put(200,100){\circle*{5}}
\put(200,100){\vector(2,1){40}}
\put(195,115){$X$}
\put(302,152){$C_1$}
\put(211,196){$C_2$}
\put(90,158){$C_3$}
\put(90,37){$C_4$}
\put(200,-10){$C_5$}
\put(302,37){$C_6$}
\put(100,50){\circle*{3}}
\put(100,150){\circle*{3}}
\put(300,50){\circle*{3}}
\put(300,150){\circle*{3}}
\put(200,0){\circle*{3}}
\put(200,200){\circle*{3}}

\multiput(100,50)(200,0){2}{\vector(0,1){8}}
\multiput(100,60)(200,0){2}{\vector(0,1){18}}
\multiput(100,80)(200,0){2}{\vector(0,1){38}}
\multiput(100,120)(200,0){2}{\vector(0,1){18}}
\multiput(100,140)(200,0){2}{\vector(0,1){8}}

\multiput(100,150)(100,-150){2}{\vector(2,1){8}}
\multiput(110,155)(100,-150){2}{\vector(2,1){18}}
\multiput(130,165)(100,-150){2}{\vector(2,1){38}}
\multiput(170,185)(100,-150){2}{\vector(2,1){18}}
\multiput(190,195)(100,-150){2}{\vector(2,1){8}}

\multiput(100,50)(100,150){2}{\vector(2,-1){8}}
\multiput(110,45)(100,150){2}{\vector(2,-1){18}}
\multiput(130,35)(100,150){2}{\vector(2,-1){38}}
\multiput(170,15)(100,150){2}{\vector(2,-1){18}}
\multiput(190,5)(100,150){2}{\vector(2,-1){8}}

\qbezier(100,50)(175,150)(300,150)
\put(148,100){\vector(2,3){10}}
\put(187,126){\vector(2,1){25}}
\put(219,138){\vector(2,1){15}}

\qbezier(100,50)(215,70)(300,150)
\put(188,76){\vector(4,1){15}}
\put(208,85){\vector(3,1){30}}
\put(235,100){\vector(2,1){15}}

\qbezier(100,50)(130,130)(170,150)
\qbezier(170,150)(210,170)(300,150)
\put(170,150){\vector(2,1){30}}

\qbezier(100,50)(205,27.5)(230,40)
\qbezier(230,40)(260,55)(300,150)
\put(230,40){\vector(2,1){30}}

\qbezier(100,50)(175,50)(215,70)
\qbezier(215,70)(255,90)(300,150)
\put(215,70){\vector(2,1){30}}
\put(122,100){\vector(1,3){8}}
\put(236,160){\vector(4,1){15}}
\put(180,58){\vector(1,0){15}}
\put(170,38){\vector(3,-1){15}}
\put(257,100){\vector(2,1){15}}
\put(278,100){\vector(1,1){10}}

\qbezier(100,50)(105,120)(120,140)
\put(120,140){\vector(1,3){4}}
\qbezier(120,140)(140,160)(160,170)
\put(160,170){\vector(2,1){30}}
\qbezier(160,170)(180,180)(200,180)
\put(200,180){\vector(3,1){10}}
\qbezier(200,180)(220,180)(300,150)
%\qbezier(120,140)(130,150)(160,160)
%\qbezier(160,160)(
%\qbezier(100,50)(100,130)(160,170)

\end{picture}
     \caption{The $2$-cell corresponding to $(X, C_1)$}\label{fig:2-cell}
\end{figure}
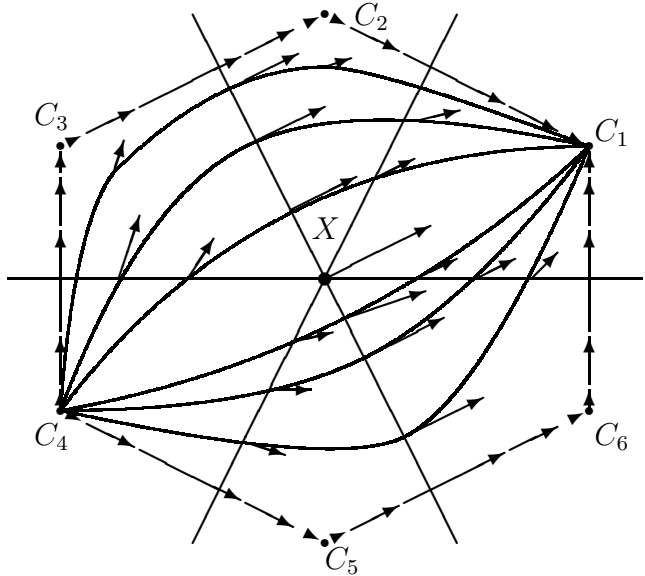

\end{itemize}
\end{Example}

\subsection{The Deligne groupoid and its representation}

In \S\ref{sec:twist} we will discuss the chain complex with 
coefficients in a local system. For the purposes, the structure of 
the fundamental group $\pi_1(\sfM(\calA))$ is particularly important. 
The concept of ``the Deligne groupoid'' $\rmGal(\calA)$ 
for a real arrangement 
$\calA$, introduced by P. Deligne \cite{del-simp} 
see also \cite{par-fundam}, and its representations 
are good tools for 
extracting information about the fundamental groups and 
local systems. 

A sequence $C_0, C_1, \ldots, C_n$ of chambers is a 
{\it gallery $G$} of {\it length n} ({\it from $C_0$ to $C_n$}) 
if $C_i$ and $C_{i+1}$ are adjacent for $i=0, 1, \ldots, n-1$. 
Any continuous path in 
$U=V_\bbR\setminus\bigcup_{X\in\calL, \codim X\ge 2}X$ 
which is transverse to any codimension one faces 
determines a gallery and every 
gallery arises in this way. Any two chambers can be connected by galleries. 
The distance between two chambers $C$ and $C'$ is the 
length of a shortest gallery connecting them; equivalently, 
it is the number of hyperplanes separating $C$ and $C'$. 
A gallery is said to be geodesic, or minimal, 
if its length is equal to the distance 
between the initial and terminal chambers. 

\begin{Def}\normalfont
\label{def:gal}
{\normalfont\bf \cite{del-simp}}
\begin{itemize}
\item[{\normalfont (1)}]
Let $G=(C_0, C_1, \ldots, C_m)$ and 
$G'=(C'_0, C'_1, \ldots, C'_n)$ be two galleries. If $C_m=C'_0$, 
define the composition of $G$ and $G'$ by 
$G G':=(C_0, \ldots, C_m=C'_0, C'_1, \ldots, C'_n)$. 

\item[{\normalfont (2)}]
Two galleries $G$ and $G'$ which have the same initial and 
terminal chambers are called {\it equivalent}, denoted by $G\sim G'$, 
if there exists a sequence of galleries $G=G_0, G_1, \ldots, G_N=G'$ 
such that for each $i=0, \ldots, N-1$, $G_i$ and $G_{i+1}$ have 
expressions
\begin{eqnarray*}
G_i&=&E_1 F E_2\\
G_{i+1}&=&E_1 F' E_2, 
\end{eqnarray*}
where $F$ and $F'$ are geodesic galleries connecting the same 
initial and terminal chambers. 

\item[{\normalfont (3)}]
$\rmGal^+(\calA)$ is defined to be the category whose objects are chambers 
$\sfch(\calA)$ and morphisms are 
$$
\Hom_{\rmGal^+}(C, C')=
\{\mbox{\normalfont Galleries from $C$ to $C'$}\}/\sim. 
$$
Since a composition of galleries is compatible with $\sim$, 
compositions of $\Hom_{\rmGal^+}$ is well-defined. 

\item[{\normalfont (4)}]
The {\it Deligne groupoid} is a category $\rmGal(\calA)$ with a functor 
$Q:\rmGal^+(\calA)\rightarrow \rmGal(\calA)$ such that 
\begin{itemize}
\item $Q(s)\in\Hom_{\rmGal}$ is an isomorphism for every $s\in\Hom_{\rmGal^+}$. 
\item Any functor $\Psi:\rmGal^+\rightarrow \calC$ such that $\Psi(s)$ 
is an isomorphism for all $s\in\Hom_{\rmGal^+}$ factors uniquely through $Q$. 
\end{itemize}
\end{itemize}
\end{Def}
See \cite{par-cov} and \cite{par-fundam} more on the 
construction of $\rmGal(\calA)$. 
The Deligne groupoid $\rmGal(\calA)$ is, roughly, obtained from 
$\rmGal^+(\calA)$ by inverting all morphisms. 
If $\calA$ is a simplicial arrangement, then 
the functor 
$Q:\rmGal^+(\calA)\rightarrow \rmGal(\calA)$ is faithful 
\cite{del-simp}. However, 
it is worth noting that $Q$ is not necessarily faithful; 
moreover 
$\rmGal^+$ is not cancellative. For example, consider 
the following two galleries 
$$
G:=C_2 C_1 C_2 C_3 C_2 \ \mbox{  and  }\ 
G':=C_2 C_3 C_2 C_1 C_2
$$
in the arrangement illustrated in Figure \ref{fig:flag}. 
Obviously $G$ and $G'$ are not equivalent in $\Hom_{\rmGal^+}(C_2, C_2)$. 
But concatenations $(C_5 C_2)G$ and $(C_5 C_2)G'$ are equivalent, 
indeed, 
$$
\begin{array}{ccl}
(C_5 C_2)G&=&C_5C_2C_1C_2C_3C_2\\
&=&C_5C_4C_1C_2C_3C_2=(C_5C_4)(C_4C_1C_2C_3)(C_3C_2)\\
&=&(C_5C_4)(C_4C_7C_6C_3)(C_3C_2)=(C_5C_4C_7)(C_7C_6C_3C_2)\\
&=&(C_5C_6C_7)(C_7C_4C_1C_2)=(C_5C_6)(C_6C_7C_4C_1)(C_1C_2)\\
&=&(C_5C_6)(C_6C_3C_2C_1)(C_1C_2)=(C_5C_6C_3)(C_3C_2C_1C_2)\\
&=&(C_5C_2C_3)(C_3C_2C_1C_2)=(C_5C_2)G'. 
\end{array}
$$
Since $(C_5C_2)$ is invertible in $\rmGal(\calA)$, 
$G$ and $G'$ determine the same element in 
$\Hom_{\rmGal}(C_2, C_2)$. 

Let $C, C'\in\sfch(\calA)$. It follows from the definition that 
any geodesic connecting $C$ to $C'$ are equivalent to each other. 
So geodesics from $C$ to $C'$ determine an equivalence class. 
We denote this equivalence class by $P^+(C, C')\in\Hom_{\rmGal}(C, C')$, 
and its inverse by 
$P^-(C', C):=P^+(C, C')^{-1}\in\Hom_{\rmGal}(C', C)$. 

\begin{Example}\normalfont
In {\normalfont Figure \ref{fig:2-cell}} galleries $(C_4, C_3, C_2, C_1)$ and 
$(C_4, C_5, C_6, C_1)$ are geodesics. Hence they determine the same 
element 
$P^+(C_4, C_1)\in\Hom_{\rmGal}(C_4, C_1)$. 
\end{Example}

Let $\calG$ be a groupoid and $x$ be an object. Then 
$\Hom_\calG(x, x)$ is 
a group and called the {\it vertex group at $x$}. The vertex group 
of the Deligne groupoid $\rmGal(\calA)$ at a chamber 
is actually isomorphic to 
the fundamental group of the complexified complement $\sfM(\calA)$ 
\cite{par-cov, par-fundam}: 
$$
\Hom_{\rmGal}(C,C)\cong\pi_1(\sfM(\calA)). 
$$
Moreover we have, 

\begin{Thm}
Let $X=X(\calA)$ be the Salvetti complex as in Theorem \ref{thm:salvetti}. 
Let $\calG(X)$ be the groupoid whose objects are $0$-cells $X_0$ and 
homomorphisms are the set of homotopy equivalence 
classes of paths between two 
$0$-cells. Then $\calG(X(\calA))$ is equivalent to the Deligne 
groupoid $\rmGal(\calA)$. 
\end{Thm}

Recall that a representation $\Phi$ of a category $\calC$ is 
a functor $\Phi:\calC\rightarrow \rmVect_\bbK$ from $\calC$ to 
the category of $\bbK$-vector spaces. $\Phi$ is 
given by a vector space $\Phi_x$ for each object $x\in\calC$ and a 
linear map $\Phi_\rho:\Phi_x\rightarrow\Phi_y$ for each 
$\rho\in\Hom_\calC(x, y)$ such that 
$\Phi_{\rho_1\rho_2}=\Phi_{\rho_1}\circ\Phi_{\rho_2}$. 

Let $\calG$ be a groupoid with 
a vertex group $G_x=\Hom(x, x)$. Then 
the category of representations $\Rep(\calG)$ of 
$\calG$ is equivalent to 
the category of group representations $\Rep(G_x)$. 
Since the category of representations of the 
fundamental group of a topological space is 
equivalent to that of local systems 
over the space, we have the following result. 

\begin{Prop}
Let $\calA$ be a real arrangement. Then 
the following categories are equivalent. 
\begin{itemize}
\item 
$\Rep(\rmGal(\calA)):$ the category of representations of 
the Deligne groupoid. 
\item 
$\Rep(\pi_1(\sfM(\calA))):$ the category of representations of 
the fundamental group. 
\item 
$\Loc(\sfM(\calA)):$ the category of local systems. 
\end{itemize}
\end{Prop}
In \S \ref{sec:twist}, we will use representations of 
the Deligne groupoid instead of local systems to compute 
the boundary maps for cellular chain complexes. 
The following operator will be needed for the purpose of 
describing the cellular boundary map. 

Let $\Phi:\rmGal(\calA)\rightarrow\rmVect_\bbK$ be a representation 
of the Deligne groupoid. Given two chambers $C$ and $C'$, 
we have two extreme 
morphisms $P^\pm(C,C'):C\rightarrow C'$. Hence we have linear maps 
$$
\Phi_{P^\pm(C,C')}:\Phi(C)\longrightarrow\Phi(C'). 
$$

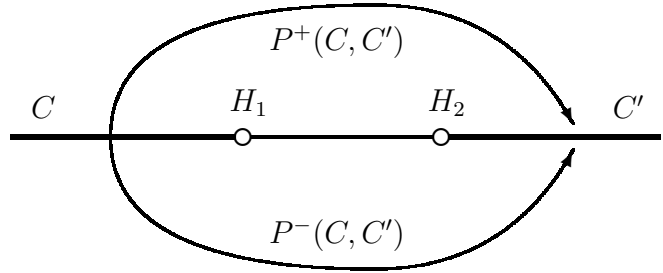
\begin{figure}[htbp]
\begin{picture}(100,100)(0,0)%(60,10)
\thicklines

\multiput(150,50)(75,0){2}{\circle{6}}
\multiput(147,49.5)(0,0.1){10}{\line(-1,0){85}}
\put(153,50){\line(1,0){69}}
\multiput(228,49.5)(0,0.1){10}{\line(1,0){80}}

\qbezier(100,50)(100,0)(200,0)
\qbezier(200,0)(250,0)(275,45)
\put(275,45){\vector(1,2){0}}

\qbezier(100,50)(100,100)(200,100)
\qbezier(200,100)(250,100)(275,55)
\put(275,55){\vector(1,-2){0}}

\put(145,60){$H_1$}
\put(220,60){$H_2$}

\put(70,57){$C$}
\put(290,57){$C'$}

\put(160,10){$P^-(C, C')$}
\put(160,83){$P^+(C, C')$}

\end{picture}
     \caption{$P^+(C, C')$ and $P^-(C, C')$}\label{fig:twohom}
\end{figure}

\begin{Def}
\label{def:diff}
\normalfont
(The {\it skein operator}). 
$$
\Delta_\Phi(C,C'):=\Phi_{P^+(C,C')}-\Phi_{P^-(C,C')}. 
$$
\end{Def}

\section{Morse theory on the complement}

\label{sec:morse}

Throughout this section, 
we investigate complex hyperplane arrangements 
which do not necessarily arise from real arrangements. 

\subsection{The Lefschetz Theorem for hyperplane complements}

\label{subsec:lef}

In this section we give a proof of the 
Lefschetz theorem for $\sfM(\calA)$ following Hamm and L\^e \cite{le-ham}. 
Although this is just a version of The Lefschetz Theorem 
for affine varieties, Morse theoretic arguments and constructions 
in this section will be needed in \S\ref{subsec:char}. 

Let $\calA=\{H_1, \ldots, H_n\}$ be an arrangement of 
hyperplanes in $\bbP_\bbC^\ell$. Let $\alpha_i$ be 
a linear form in $\bbC[z_0, z_1, \ldots, z_\ell]$ defining $H_i$ and 
$Q$ denote the product $\alpha_1\alpha_2\cdots\alpha_n$ 
of these linear forms. 
Let $V(Q)$ be the union $\bigcup_{i=1}^n H_i$ of hyperplanes and 
$\cpl(Q)=\bbP^\ell-V(Q)$ denote the complement. 
There exists an 
obvious stratification $\Sigma(\calA)$ 
of the union as follows. 
Given an intersection $X\in L(\calA)$ of some 
hyperplanes in $\calA$, define 
$$
S_X:=X-\bigcup_{H\nsupseteq X}H. 
$$
We have a partition 
$\{S_X\}_{X\in L(\calA)}$ of $\bbP^\ell$. 

\begin{Lemma}
For an arrangement $\calA$, the above stratification 
$\Sigma(\calA)=\{S_X\}_{X\in L(\calA)}$ is a good stratification 
at each point $p\in V=V(Q)$, i.e. there exist 
a neighborhood $\calU\ni p$ and a holomorphic function $h$ on $\calU$ 
with $V(h)=\calU\cap V(Q)$ satisfying the following 
Thom's condition {\normalfont ($a_h$)}: 
\begin{itemize}
\item[{\normalfont ($a_h$)}]
If $p_i$ is a sequence of points in $\calU-V(h)$ 
such that $p_i\rightarrow p\in S_X$ 
and $\rmT_{p_i}V(h-h(p_i))$ converges to some hyperplane 
$\calT$, then $\rmT_p S_X\subset \calT$. 
\end{itemize}
\end{Lemma}

The rest of this section is devoted to proving the 
following theorem (\cite{ham-aff, le-ham, dim-pap, ran-mor}). 

\begin{Thm}
\label{thm:zale}
\begin{itemize}
\item[{\normalfont (i)}] 
Let $F=V(f)\subset\bbP_\bbC^\ell$ be a hyperplane defined by a 
linear form $f$ which is transverse to each stratum. Then 
$\cpl(Q)$ has the homotopy type of a space obtained from 
$\cpl(Q)\cap F$ by attaching a certain number of $\ell$-dimensional 
cells. 
\item[{\normalfont (ii)}] Moreover the number of $\ell$-cells is the 
$\ell$-th Betti number $b_\ell(\cpl(Q))$. 
\end{itemize}
\end{Thm}
(ii) is proved in \S\ref{sec:intro}.  
The plan of the proof of (i) is to apply 
Morse theory to a function of the form 
\begin{equation}
\varphi(x)=
\left|\frac{f(x)^{\lambda_0}}
{\alpha_1^{\lambda_1}\cdots\alpha_n^{\lambda_n}}\right|^2, 
\mbox{ for }x\in \cpl(g), 
\end{equation} 
where $\lambda_1, \ldots, \lambda_n\in\bbZ_{>0}$ are 
appropriately chosen positive integers and 
$\lambda_0=\lambda_1+\cdots +\lambda_n$. Note that 
$\varphi$ is a well-defined differentiable map 
from $\cpl(Q)$ to $\bbR_{\ge 0}$ which has the bottom 
$F\cap \cpl(Q)=\varphi^{-1}(0)$. 
The reason for 
considering this function is that the critical points 
are well studied, in particular, 
critical points are known to be nondegenerate 
for generic $\lambda_1, \ldots, \lambda_n$. 
It was conjectured by Varchenko \cite{var-cri}, and 
proved by Orlik-Terao \cite{ot-num} and Silvotti \cite{sil-var}. 

\begin{Thm}
\label{thm:otvar}
Let $\calA$ be a complex essential affine arrangement in $\bbC^\ell$ 
with defining linear equations $f_1, \ldots, f_N$, 
and put 
$$
\Phi_\lambda=f_1^{\lambda_1}f_2^{\lambda_2}\cdots
f_N^{\lambda_N}
$$ 
for $\lambda=(\lambda_1, \ldots, \lambda_N)\in\bbC^N$. Then 
there exists a Zariski-closed algebraic proper subset 
$Y$ of $\bbC^N$, such that for $\lambda\in\bbC^N-Y$, 
$\Phi_\lambda$ has only finitely many critical points, 
all of which are nondegenerate and the number of critical points 
of $\Phi_\lambda$ is $|\chi(\sfM(\calA))|$. 
\end{Thm}

In our situation, since 
$$
\Phi_\lambda=\frac{f^{\lambda_0}}
{\alpha_1^{\lambda_1}\cdots\alpha_n^{\lambda_n}}=
(\alpha_1/f)^{-\lambda_1}\cdots(\alpha_n/f)^{-\lambda_n}, 
$$
there exist $\lambda_1, \ldots, \lambda_n\in\bbZ_{>0}$ 
such that $\Phi_\lambda$ has only nondegenerate critical points. 
Combining the above theorem with Proposition \ref{prop:flag}, 
the number of critical points is shown to be equal to 
the $\ell$-th Betti number $b_\ell(\sfM(Q))$ of the complement. 
From the next lemma, 
$\varphi=|\Phi_\lambda|$ also 
has only finitely many critical points all of which 
are nondegenerate critical points of Morse index $\ell$. 

\begin{Lemma}
\label{lem:morse}
Let $\frakf$ and $\frakg:U\rightarrow \bbC$ be holomorphic functions defined 
on a neighborhood $U$ of $0\in\bbC^n$. We assume $\frakf(0), \frakg(0)\neq 0$. 
\begin{itemize}
\item[{\normalfont (i)}] $0\in U$ is a critical point of $|\frakf|^2$ if 
and only if $0\in U$ is a critical point of $\frakf$. 
\item[{\normalfont (ii)}] In {\normalfont (i)}, $0\in U$ is a 
nondegenerate critical point of $|\frakf|^2$ if 
and only if $0\in U$ is a nondegenerate critical point of $\frakf$. 
\item[{\normalfont (iii)}] If $0\in U$ is a nondegenerate critical point 
of $|\frakf|^2$, then the Morse index is $n$. 
\item[{\normalfont (iv)}] If $d\frakf$ and $d\frakg$ are linearly 
independent over $\bbC$ at each point in $U$, 
then so are $d|\frakf|$ and $d|\frakg|$. 
\end{itemize}
\end{Lemma}

\noindent
{\bf Proof. }
Since $|\frakf|^2=\frakf\cdot\bar{\frakf}$, 
$\frac{\partial}{\partial z_i}|\frakf|^2=
\frac{\partial\frakf}{\partial z_i}\bar{\frakf}$ 
and 
$\frac{\partial}{\partial \bar{z_i}}|\frakf|^2=
\frac{\partial\bar{\frakf}}{\partial \bar{z_i}}\frakf$. Thus we have (i). 
Moreover the determinant of the Hessian matrix at $0\in U$ is 
\begin{eqnarray*}
\det
\left(
\begin{array}{cc}
\frac{\partial^2 |\frakf|^2}{\partial z_i\partial z_j}
&
\frac{\partial^2 |\frakf|^2}{\partial z_i\partial \bar{z_j}}
\\
\frac{\partial^2 |\frakf|^2}{\partial \bar{z_i} \partial z_j}
&
\frac{\partial^2 |\frakf|^2}{\partial \bar{z_i} \partial \bar{z_j}}
\end{array}
\right)
&=&
\det
\left(
\begin{array}{cc}
\frac{\partial^2 \frakf}{\partial z_i\partial z_j}\bar{\frakf}
&
\frac{\partial \frakf}{\partial z_i}
\frac{\partial \bar{\frakf}}{\partial \bar{z_j}}
\\
\frac{\partial \bar{\frakf}}{\partial \bar{z_i}}
\frac{\partial \frakf}{\partial z_j}
&
\frac{\partial^2 \bar{\frakf}}{\partial \bar{z_i} \partial \bar{z_j}}\frakf
\end{array}
\right)
\\
&=&
|\frakf|^2 
\left|\det
\left(
\frac{\partial^2 \frakf}{\partial z_i\partial z_j}
\right)
\right|^2. 
\end{eqnarray*}
(Here we use 
$\frac{\partial \frakf}{\partial z_i}=
\frac{\partial\bar{\frakf}}{\partial\bar{z_i}}=0$. ) 
This proves (ii). 

After a linear change of coordinates, we may assume $\frakf$ is 
expressed as 
$$
\frakf(z_1, \ldots, z_n)=c(1+\sum_{i=1}^n z_i^2+O(3))
$$
with $c\neq 0$. Writing $z_i=x_i+\ii y_i$, we have 
\begin{eqnarray*}
\frakf\bar{\frakf}&=&|c|^2(1+\sum_{i=1}^n(z_i^2+\bar{z_i}^2) +O(3))\\
&=&|c|^2(1+2\sum_{i=1}^n(x_i^2-y_i^2) +O(3)). 
\end{eqnarray*}
Consequently the Morse index of $|\frakf|^2$ at $0$ is equal to $n$. 

(iv) is clear from 
$d(|\frakf|^2)=(d\frakf)\bar{\frakf}+\frakf(d\bar{\frakf})$. 
\owari

%%%%%%%%%%%%%%%%%%%%%%%%%%%%%%%%%%%%%%%%%%%%%%%%%%%%%%%%%%%%%%%%%%%%%
%%%%%%%%%%%%%%%%%%%%%%%%%%%%%%%%%%%%%%%%%%%%%%%%%%%%%%%%%%%%%%%%%%
%%%%%%%%%%%%%%%%%%%%%%%%%%%%%%%%%%%%%%%%%%%%%%%%%%%%%%%%%%%%%%%%%%

Unfortunately, our Morse function 
$\varphi=|\Phi_\lambda|$ is not a proper function. 
Hence it is necessary to study the Morse theory for a 
nonproper Morse function. This difficulty is directly 
related to the fact that $\varphi$ 
has points of indeterminacy: $V(f, Q)=\{f=Q=0\}=V(Q)\cap F$. 
To deal with the difficulty, we have to remove a neighborhood 
of $V(Q)\cap F$ for separating the zero loci and the poles of $\varphi$. 
The idea is to measure the distance from $V_1:=V(f, Q)$. 
Suppose $p=[z_0:\cdots :z_\ell]\in\bbP^\ell-V_1$, and 
define $h_{V_1}(p)$ as follows 
$$
h_{V_1}(p)=\frac{|z_0|^{\lambda_0}+|z_1|^{\lambda_0}+\cdots +|z_\ell|^{\lambda_0}}
{|f|^{\lambda_0}+|\alpha_1^{\lambda_1}\cdots\alpha_n^{\lambda_n}|}.  
$$
Then 
$h_{V_1}:(\bbP^\ell-V_1)\rightarrow\bbR_{\ge 0}$ is a well-defined 
map. 
%%Moreover for a converging sequence $p_i\rightarrow p\in V_1$, 
%%we have $h_{V_1}(p_i)\rightarrow\infty$. 

\begin{Lemma}
\label{lem:diffeo}
Let 
$M^{\leq t}:=\{p\in\bbP^\ell-V_1; h_{V_1}(p)\leq t\}$. 
For sufficiently large $t\gg 0$, 
$h_{V_1}^{-1}(t)=\partial M^{\leq t}$ is 
transverse to each stratum $S\in\Sigma$ and to $F$. Moreover 
$\left(\bbP^\ell-V_1, (\bbP^\ell-V_1)\cap\Sigma\right)$ 
is diffeomorphic to $\left(M^{< t}, M^{< t}\cap \Sigma\right)$ 
as stratified manifolds. 

\begin{figure}[htbp]
\begin{picture}(100,150)(0,0)
\thicklines

\multiput(70,30)(20,0){13}{\line(1,0){15}}

\put(300,30){\circle*{4}}
\put(300,30){\line(1,1){10}}
\put(310,45){$F=\{f=0\}$}

\put(105,30){\line(3,4){95}}
\put(105,30){\line(-3,-4){20}}
\put(105,30){\circle*6{6}}

\put(252,30){\line(-4,5){100}}
\put(252,30){\line(4,-5){20}}
\put(252,30){\circle*6{6}}

\put(150,45){\vector(-3,-1){37}}
\put(200,45){\vector(4,-1){45}}
\put(147,48){$V_1=F\cap V$}

\put(188,110){\circle*{4}}
\put(188,110){\line(3,2){17}}
\put(188,140){\circle*{4}}
\put(188,140){\line(3,-2){17}}
\put(209,120){$V=\{Q=0\}$}

\put(105,30){\circle{36}}
\put(252,30){\circle{44}}

\put(150,2){$\partial M^{\leq t}=h_{V_1}^{-1}(t)$}
\put(156,13){\vector(-4,1){35}}
\put(205,13){\vector(3,1){28}}

\end{picture}
     \caption{}
\label{fig:zl}
\end{figure}
\end{Lemma}

\noindent
{\bf Proof.} 
We first observe that $h_{V_1}$ is defined on 
$\bbP^\ell-V_1$ with values in $\bbR_{>0}$. 
It is clear that for a sequence $\{p_i\}\subset\bbP^\ell-V_1$ 
converging to a point $p\in V_1$, we have 
$h_{V_1}(p_i)\rightarrow\infty$ 
as $i\rightarrow\infty$. Recall that any real polynomial function 
on a semi-algebraic set can have at most a finite number of 
critical values (see Milnor \cite[Cor. 2.8]{mil-sing}). 
Since any restriction $h_{V_1}|_S$ to a stratum $S\in\Sigma$ 
has only finitely many critical values, we may choose $t$ to be 
larger than any critical value. 
Suppose $t_1$ and $t_2$ ($t_1<t_2$) are sufficiently 
large. Then there exists an open neighborhood $U$ of 
the compact set $M^{[t_1, t_2]}:=\{p; t_1\leq h_{V_1}(p)\leq t_2\}$ 
such that the restriction of $h_{V_1}$ to $U\cap S$ has no 
critical points for any stratum $S\in\Sigma$. 

The gradient vector field $-\grad h_{V_1}$ does not 
preserve the stratification in general. 
We modify $-\grad h_{V_1}$ so that it preserves 
the stratification. 
Let 
$p\in U$
and $S$ denote the stratum which contains $p$. Since $p$ is 
not a critical point of $h_{V_1}|_S$, there exists 
a tangent vector $v\in\rmT_pS$ such that $v\cdot\varphi<0$. 
On a small neighborhood $U_p$ of $p$ in $U$, not meeting any smaller 
stratum than $S$,  
we have a vector field $\tilde{v}$ such that 
\begin{itemize}
\item[(i)] $\tilde{v}$ is tangent to each stratum $\Sigma\cap U_p$, 
\item[(ii)] $\tilde{v}\cdot\varphi<0$. 
\end{itemize}
Note that since our strata are linear, we can take $\tilde{v}$ 
as a constant vector field in a certain open set. 

Using a partition of unity, we have a vector field $\tilde{v}$ on 
$U$ satisfying the conditions (i) and (ii) above. 
Then $\tilde{v}/||\tilde{v}||$ defines a deformation retract of 
$M^{<t_2}$ onto $M^{<t_1}$, which preserves 
the structure of stratification. 
\owari

Now we consider the function $\varphi$ on $M^{\leq t}$ and the restriction 
to its boundary $\partial M^{\leq t}=h_{V_1}^{-1}(t)$. 
The following lemma plays a key role in the arguments below. 
The assumption that $F$ is generic is used in the proof. 

\begin{Lemma}
\label{lem:nocrit}
For sufficiently large 
$t\gg 0$, $\varphi|_{\partial M^{\leq t}\setminus (V(Q)\cup F)}$ 
has no critical points. 
\end{Lemma}

\noindent
{\bf Proof. } Let $p\in V_1=V\cap F$ and 
$S_X\in\Sigma$ be the stratum containing $p$. Note that 
$X\in L(\calA)$ is the smallest intersection containing $p$, 
and it is, by definition, transverse to $F$. We have 
coordinates $(z_1, \ldots, z_\ell)$ in a neighborhood 
$U$ of $p$ with the origin at $p$. The transversality 
of $F$ to $\Sigma$ allows us to assume that 
\begin{itemize}
\item[(1)] 
$F\cap U$ is defined by a linear form $z_\ell=0$. 
\item[(2)] 
$X\cap U$ is defined by $\{z_1=z_2=\cdots =z_m=0\}$ with $1\leq m<\ell$. 
\item[(3)] 
Let $H_1, \ldots, H_k$ be the set of all 
hyperplanes in $\calA$ which contains $p$. Each $H_i\cap U$ 
is defined by a linear 
form of the form $a_1z_1+\cdots + a_mz_m$. 
\end{itemize}
For simplicity, set $g_1=\alpha_1^{\lambda_1}\cdots\alpha_k^{\lambda_k}$ and 
$g_2=z_\ell^{\lambda_0}$. 
The assumptions imply that $d|g_1|$ and $d|g_2|$ are linearly 
independent at each point of $U-(V(Q)\cup F)$. 
Now $h_{V_1}$ and $\varphi$ are expressed as 

\begin{eqnarray*}
h_{V_1}&=&\frac{1+|z_1|^{\lambda_0}+\cdots +|z_\ell|^{\lambda_0}}
{|g_1|+|g_2|}\\
\varphi&=&\frac{|g_2|}{|g_1|}. 
\end{eqnarray*}
Now we prove that there exists a neighborhood $U'$ of $p$ such that 
$dh_{V_1}$ and $d\varphi$ are linearly independent at each point 
of $U'-(V(Q)\cup F)$. 
\begin{eqnarray}
\label{eq:dh}
d\log h_{V_1}&=&-\frac{d|g_1|+d|g_2|}{|g_1|+|g_2|}
+d\log(1+|z_1|^{\lambda_0}\cdots +|z_\ell|^{\lambda_0})\\
d\log\varphi&=&-\frac{d|g_2|}{|g_2|}+\frac{d|g_1|}{|g_1|}. 
\end{eqnarray}
If $U'$ is sufficiently small, then the last term of (\ref{eq:dh}) 
is sufficiently small. Compare the sign of coefficients of 
$d|g_1|$ and $d|g_2|$, we conclude that $dh_{V_1}$ and $d\varphi$ 
are linearly independent. 
Thus for any point $p\in F\cap V(Q)$, there exists a neighborhood 
$U_p$ in $\bbP_\bbC^\ell$ such that $d\varphi$ and $dh_{V_1}$ are 
linearly independent at each point of $U_p-(V(Q)\cup F)$. 
We choose finitely many points $p_1, \ldots, p_N\in V(Q)\cap F$ with 
$V(Q)\cap F\subset \bigcup_{i=1}^N U_{p_i}$ and set 
$t_0:=\sup\{h_{V_1}(p);\ p\in \bbP_\bbC^\ell-\bigcup_{i=1}^N U_{p_i}\}$. 
Then for $t>t_0$, 
$\varphi|_{\partial M^{\leq t}\setminus (V(Q)\cup F)}$ 
has no critical points. 
\owari 

The transversality $F\pitchfork\partial M^{\leq t}$ is also 
shown by observing (\ref{eq:dh}). Hence a small tubular neighborhood 
of $F\cap\partial M^{\leq t}$ in $\partial M^{\leq t}$ is 
diffeomorphic to a disk bundle over $F\cap\partial M^{\leq t}$. 
Recall that the normal bundle of $F$ in $\bbP_\bbC^\ell$ is 
$\calN_{F/\bbP_\bbC^\ell}\cong \calO_F(1)$ and is trivial 
on an affine open set $F-V(Q)$. Thus we have: 
\begin{Lemma}
\label{lem:tub}
There exists a small neighborhood $\calT$ of $F$ in $\bbP_\bbC^\ell$ 
such that 
\begin{equation}
\label{eq:tub}
\calT\cap M^{\leq t}\cong (F\cap M^{\leq t})\times \rmD^2\ 
\mbox{\normalfont  (Diffeomorphic),}
\end{equation}
where $\rmD^2=\{(x,y)\in\bbR^2; x^2+y^2<1\}$ is the unit disk. 
\end{Lemma}

\begin{figure}[htbp]
\begin{picture}(100,130)(0,0)
\thicklines

\multiput(0,0)(15,0){12}{\line(1,0){10}}
\put(0,3){$F$}
%\qbezier(149,-8)(150,100)(170,130)
\qbezier(149,-8)(170,130)(170,130)
\put(175,125){$V=V(g)$}
\put(150,0){\circle*{6}}
\put(154,4){$F\cap V$}
\qbezier(170,100)(100,100)(70,50)
\qbezier(70,50)(52,20)(50,-8)
\put(175,95){$\partial N$}
\put(50,109){$X=-\grad \varphi$}

\put(150,90){\vector(-1,1){15}}
\put(130,90){\vector(-2,1){18}}
\put(110,80){\vector(-1,0){18}}
\put(90,70){\vector(-1,-1){18}}
\put(70,60){\vector(0,-1){20}}
\put(60,40){\vector(1,-2){10}}
\put(50,20){\vector(1,-1){14}}

\multiput(250,0)(15,0){12}{\line(1,0){10}}
\put(250,3){$F$}
\qbezier(369,-8)(390,130)(390,130)
\put(395,125){$V$}
\put(370,0){\circle*{6}}
\put(374,4){$F\cap V$}
\qbezier(390,100)(320,100)(290,50)
\qbezier(290,50)(272,20)(270,-8)
\put(395,95){$\partial N$}
\put(270,90){$\calX$}

\put(375,100){\vector(-1,0){20}}
\put(340,92){\vector(-4,-1){20}}
\put(323,83){\vector(-3,-1){20}}
\put(305,70){\vector(-2,-1){20}}
\put(290,50){\vector(-2,-3){15}}
\put(275,20){\vector(-1,-2){10}}

\end{picture}
     \caption{Modified vector field}\label{fig:modify}
\end{figure}
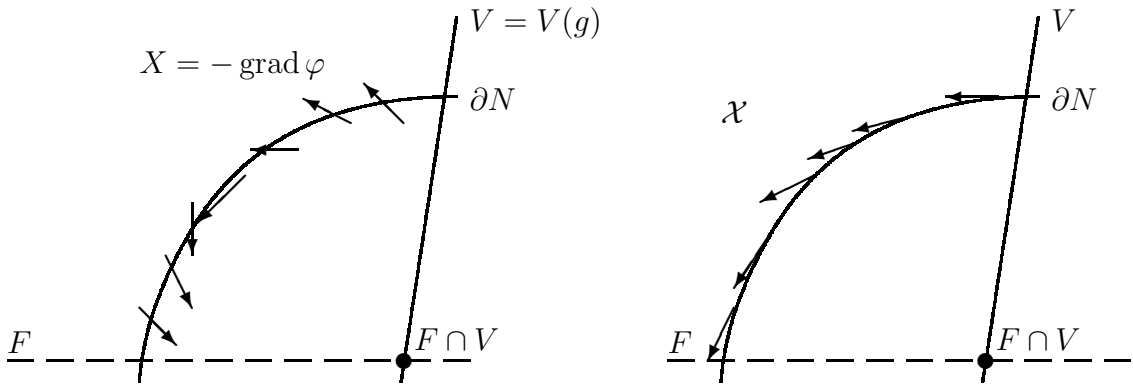

Now we return to the proof of Lefschetz's Theorem \ref{thm:zale}. 
Fix a sufficiently large $t\gg 0$ and set $N=M^{\leq t}\cap \cpl(Q)$. 
Consider 
the gradient vector field $X=-\grad\varphi$. 
$X$ is not tangent to the boundary $\partial N$ in general, 
so we comb the vector field in order to make it neat as in 
the proof of Lemma \ref{lem:diffeo}. 
Lemma \ref{lem:nocrit} means that $d\phi$ is not orthogonal to 
$\partial N$. Thus by an argument 
similar to that in the proof of Lemma \ref{lem:diffeo}, 
modifying $X$ around $\partial N$, we have 
a vector field $\calX$ on $N$ which is tangent to $\partial N$. 
Moreover, there exists a vector field $\calX$ 
which satisfies the following conditions 
(Fig. \ref{fig:modify}):
\begin{itemize}
\item[(a)] There exists a neighborhood $U$ of $\partial N\cup (N\cap F)$ 
such that $\calX=X=-\grad \varphi$ outside $U$ and 
$\rmCrit(\varphi)\cap \overline{U} =\emptyset$, where 
$\rmCrit(\varphi)$ is the set of critical points of $\varphi$ 
in $\cpl(Q)-F$. 
\item[(b)] $\calX$ is tangent to  $\partial N$. 
\item[(c)] $\calX\varphi<0$ on $N-(F\cup \rmCrit(\varphi))$. 
\item[(d)] Under the diffeomorphism (\ref{eq:tub}), the vector field $\calX$ 
coincides with the negative vertical Euler vector field 
$-x\frac{\partial}{\partial x}-y\frac{\partial}{\partial y}$. 
\end{itemize}

We now complete the proof of Theorem \ref{thm:zale}. 
We consider $\calX$ as the negative gradient vector field of 
a Morse function on $N$ and 
define 
$$
N^{\leq s}:=\{p\in N; \varphi(p)\leq s\}, 
$$
for $s>0$. If there is no critical value in the interval $[s_1, s_2]$, 
then the vector field $\calX$ induces a retraction 
$N^{\leq s_2}\stackrel{\cong}{\longrightarrow}N^{\leq s_1}$. 
If there is only one critical point within the interval 
$[s_1, s_2]$, then the homotopy type of 
$N^{\leq s_2}$ is obtained from that of $N^{\leq s_1}$ 
by attaching an $\ell$-cell. 
This completes the proof of Lefschetz's hyperplane section theorem. 
\owari

\subsection{Stable and unstable manifolds}

Consider the flow of $\calX$: 

\begin{eqnarray*}
&&\phi_t:N\longrightarrow N,\ t\in\bbR\\
&&\frac{\partial}{\partial t}\phi_t(x)=\calX_{\phi_t(x)},\ \phi_0={\mathrm{id}}_N.
\end{eqnarray*}
If $p\in N-F$ is a critical point of $\varphi$, we define the 
stable manifold $W_p^s$ and unstable manifold $W_p^u$ as 
\begin{eqnarray*}
W_p^s&=&\{ x\in N;\ \lim_{t\rightarrow  \infty}\phi_t(x)=p\}\\
W_p^u&=&\{ x\in N;\ \lim_{t\rightarrow -\infty}\phi_t(x)=p\}. 
\end{eqnarray*}
These are $\ell$-dimensional (over $\bbR$) submanifolds in $N$. 
Recall that 
the vector field $\calX$ is said to satisfy the {\it Morse-Smale condition} 
if the stable and unstable manifolds intersect transversely. 
But in our case, since any critical point in $N-F$ has the 
middle index $\ell$, 
it seems reasonable to define as follows. 

\begin{Def}
\normalfont
Let $\calX$ be a vector field on $N$ as in the previous section, 
it is said to satisfy the {\it Morse-Smale condition} if 
there does not exist a flow line connecting distinct points 
in $\rmCrit(\varphi)$. In other words, there does not exist 
$x\in N-(F\cup\rmCrit(\varphi))$ such that both 
$\lim_{t\rightarrow\infty}\phi_t(x)$ and 
$\lim_{t\rightarrow -\infty}\phi_t(x)$ are 
contained in $\rmCrit(\varphi)$. 
\end{Def}

Recall that for a critical point $p\in\rmCrit(\varphi)$, 
$N^{\leq \varphi(p)+\epsilon}$ is homotopy equivalent to 
$N^{\leq \varphi(p)-\epsilon}\cup W_p^u$ (see 
\cite[Thm. 3.2]{mil-mor}) and that 
unstable manifolds are 
preserved by the action of $\phi_t$. 
Hence under the Morse-Smale condition, the boundary of an 
unstable manifold $W_p^u$ should be attached to 
$N^0\subset F\cap \cpl(Q)$. 
Thus we have: 

\begin{Thm}
\label{thm:smale}
If $\calX$ satisfies the Morse-Smale condition, then $\cpl(Q)$ is 
homotopy equivalent to 
$$
(F_\bbC\cap \cpl(Q))\cup\bigcup_{p\in\rmCrit(\varphi)}W_p^u. 
$$
\end{Thm}

\begin{figure}[htbp]
\begin{picture}(100,180)(0,0)
\thicklines

\put(0,20){\line(0,1){100}}
\put(10,0){\line(0,1){100}}
\put(190,145){\line(0,-1){100}}
\put(180,165){\line(0,-1){20}}
\put(10,0){\line(4,1){180}}
\qbezier(0,20)(4,21)(8,22)
\qbezier(0,120)(40,130)(45,120)
\qbezier(45,120)(50,110)(10,100)
\qbezier(180,165)(120,150)(125,140)
\qbezier(125,140)(130,130)(190,145)
\qbezier(45,120)(60,90)(100,100)
\put(100,100){\circle*{3}}
\put(90,105){$p_2$}
\qbezier(100,100)(120,105)(124.5,142)
\put(150,157){\circle*{3}}
\put(166,140){\circle*{3}}
\qbezier(150,157)(150,177)(158,168.5)
\qbezier(158,168.5)(166,160)(166,140)

\put(158,168.5){\circle*{3}}
\put(160,175){$p_1$}
%\put(160,125){$W_{p_1}^u$}
\put(125,167){$W_{p_1}^u$}

\put(90,50){$N^{\leq s_2}$}

\put(195,50){$\approx$}

\put(220,0){\line(0,1){40}}
\put(220,0){\line(4,1){180}}
\put(210,20){\line(0,1){40}}
\qbezier(210,20)(214,21)(218,22)
\put(220,40){\line(4,1){180}}
\put(210,60){\line(4,1){180}}
\put(400,45){\line(0,1){40}}
\put(390,105){\line(0,-1){20}}

\put(290,80){\circle*{3}}
\qbezier(290,80)(290,110)(295,100)
\put(295,100){\circle*{3}}
\qbezier(295,100)(300,90)(300,60)
\put(300,60){\circle*{3}}

\put(350,157){\line(0,-1){62}}
\qbezier(350,157)(350,177)(358,168.5)
\put(358,168.5){\circle*{3}}
\qbezier(358,168.5)(366,160)(366,140)
\put(366,140){\line(0,-1){64}}

%%%%%%%%%%%%%%
%\put(350,95){\circle*{3}}
%\qbezier(350,95)(350,180)(355,170)
%\put(355,170){\circle*{3}}
%\qbezier(355,170)(360,160)(360,75)
%\put(360,75){\circle*{3}}
%%%%%%%%%%%%%%

\put(280,115){$W_{p_2}^u$}
\put(325,167){$W_{p_1}^u$}

\put(300,35){$N^{\leq s_1}$}

\end{picture}
     \caption{Unstable manifolds}\label{fig:flow}
\end{figure}
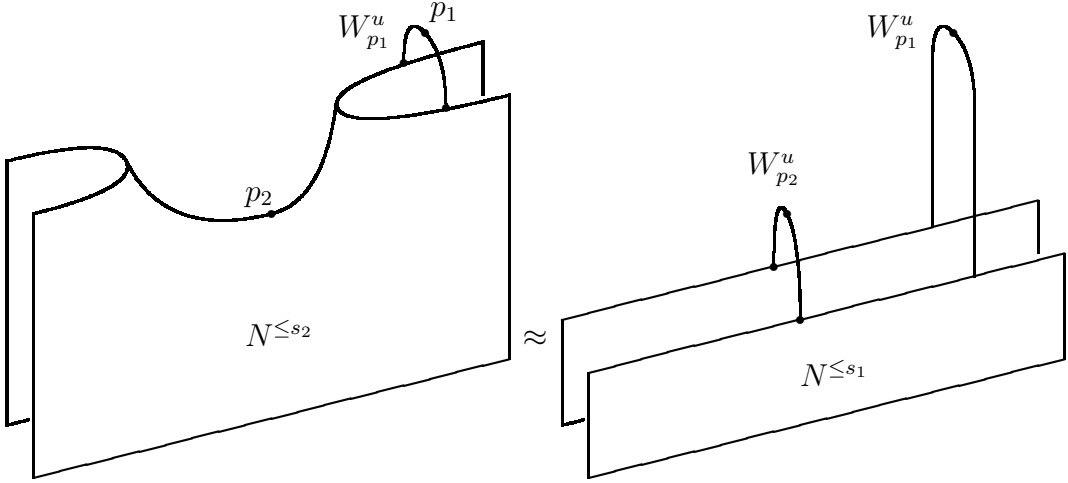

\begin{Thm}
\label{thm:diffeo}
$\cpl(Q)-\bigcup_{p\in\rmCrit(\varphi)}W_p^s$ 
is diffeomorphic to $(\cpl(Q)\cap F_\bbC)\times \rmD^2$. 
\end{Thm}

\noindent
{\bf Proof. } 
Let $\calT$ be a tubular neighborhood of $F_\bbC\cap \cpl(Q)$ in $\cpl(Q)$ 
as in Lemma \ref{lem:tub}. 
Since $\calX$ is a complete vector field, 
$$
\begin{array}{ccc}
(\partial\calT\cap N)\times \bbR&\longrightarrow &
N-\left(F_\bbC\cup\bigcup_{p\in\rmCrit(\phi)}W_p^s\right)\\
&&\\
(q, t)&\longmapsto&\phi_t(q)
\end{array}
$$
defines a diffeomorphism. It is also diffeomorphic to 
$(F_\bbC\cap N)\times (\rmD^2-\{(0,0)\})$. 
The condition (d) on $\calX$ allows us to complete the proof. 
\owari

\subsection{Homotopy types of the unstable cells}

\label{subsec:char}

Next we characterize the homotopy types of unstable 
manifolds $\{W_p^u;\ p\in\rmCrit(\varphi)\}$. 
The unstable manifold $W_p^u$ corresponding to 
$p\in\rmCrit(\varphi)$ can be considered as an attached cell. 
There exists a continuous map 
$\sigma_p:(\rmD^\ell, \partial\rmD^\ell)\rightarrow (\cpl(Q), F\cap \cpl(Q))$ 
such that $\sigma_p(0)=p$ and $\sigma_p$ induces a 
diffeomorphism of $\rmint(\rmD^\ell)$ to $W_p^u$. 
We now assume that our manifolds are oriented. 
Observe that $\sigma_p$ satisfies the following properties: 
\begin{itemize}
\item[(i)] $\sigma_p(0)=p$ and $\sigma_p(\rmD^\ell)\cap W_p^s=\{p\}$. 
\item[(ii)] $\sigma_p(\rmD^\ell)$ intersects $W_p^s$ at $p$ transversally 
and positively. 
\item[(iii)] $\sigma_p(\partial\rmD^\ell)\subset F\cap \cpl(Q)$. 
\item[(iv)] If $q\in\rmCrit(\varphi)\setminus\{p\}$ is another 
critical point, then $\sigma_p(\rmD^\ell)$ does not intersect $W_q^s$. 
\end{itemize}
Note that (iv) is equivalent to the Morse-Smale condition 
($W_p^u\cap W_q^s=\emptyset$). 
Let us call these properties 
``set-theoretical duality'' between 
cells $\{\sigma_p\}_{p\in\rmCrit(\varphi)}$ and 
stable manifolds $\{W_p^s\}_{p\in\rmCrit(\varphi)}$. 
The main result of this section is to characterize 
the homotopy type of the map 
$\sigma_p:(\rmD^\ell, \partial\rmD^\ell)\rightarrow (N, N\cap F)$ 
by set-theoretical duality for stable manifolds. 

\begin{Thm}
\label{thm:char}
Suppose that a continuous map 
$\sigma'_p:(\rmD^\ell, \partial\rmD^\ell)\rightarrow (N, N\cap F)$ 
is differentiable in a neighborhood of $0\in\rmD^\ell$ and 
satisfies conditions {\normalfont (i)} through {\normalfont (iv)} above. 
Then 
$\partial \sigma_p$ and $\partial \sigma'_p: \partial\rmD^\ell\rightarrow N\cap F$ 
are homotopic. In particular, 
$$
\sfM(\calA) \mbox{\normalfont \ and }  (\sfM(\calA)\cap F_\bbC)
\cup_{(\partial\sigma'_p)}
\left(
\coprod_{p\in\rmCrit(\phi)}\rmD^\ell
\right)
$$
are homotopy equivalent. 
\end{Thm}

\noindent
{\bf Proof. } 
The idea of the proof is simple: flowing $\sigma'_p$ via the gradient 
flow $\phi_t$, then $\phi_t\circ\sigma'_p$ 
converges to $\sigma_p$ as $t\rightarrow\infty$. 

From (i), we have $\sigma'_p(0)=\sigma_p(0)=p$ and 
the image $\sigma'_p(\rmD^\ell)$ is transverse to $W_p^s$. 
Note that $\rmT_p\bbP_\bbC^\ell=\rmT_p W_p^u\oplus\rmT_p W_p^s$. 
And the projection $\rmT_p\bbP^\ell\rightarrow \rmT_p W_p^u$ induces 
an orientation preserving isomorphism 
$\rmT_p\sigma'_p(\rmD^\ell)\cong \rmT_p W_p^u$. 
By modifying $\sigma'_p$ up to homotopy, we have $\sigma''_p$ 
satisfying (i) $\cdots$ (iv) and the following properties: 
$$
\sigma''_p(x)=
\left\{
\begin{array}{ll}
\sigma_p(x)   &\mbox{ if } ||x||<\epsilon\\
\sigma'_p(x)  &\mbox{ if } 2\epsilon\leq ||x||\leq 1, 
\end{array}
\right.
$$
where 
$||x||^2=x_1^2+\cdots +x_\ell^2$ 
and 
$\epsilon$ is a sufficiently small positive number. 
Take a tubular 
neighborhood $\calT$ of $F\cap N$ such that 
$\calT=(F\cap N)\times \rmD^2$ as in  Lemma \ref{lem:tub}. 
Denote by $\pi: \calT=(F\cap N)\times \rmD^2\rightarrow F\cap N$ the 
projection. 
Consider $\phi_t\circ \sigma''_p$.  
If $t\gg 0$ is sufficiently large 
then we may assume that $\phi_t\circ \sigma''_p(x)\in\calT$ 
for $\epsilon\leq ||x||\leq 1$. 
By definition, $(\pi\circ\phi_t\circ \sigma''_p)|_{||x||=\epsilon}$ 
is equal to $\partial \sigma_p=\sigma_p|_{\partial \rmD^\ell}$ as maps 
$\rmS^\ell\rightarrow F\cap N$, more precisely, for $x\in\partial\rmD^\ell$, 
$(\pi\circ\phi_t\circ \sigma''_p)(\epsilon x)=\partial\sigma_p(x)$. 
Since 
$(\pi\circ\phi_t\circ \sigma''_p)|_{||x||=1}=\partial \sigma'_p$, 
$h_r(x):=\pi\circ\phi_t\circ \sigma''_p(r\cdot x)$ for 
$\epsilon\leq r\leq 1$ 
defines a homotopy between 
$\partial \sigma_p$ and $\partial \sigma'_p$. 
%Then $\pi\circ\phi_t\circ \sigma''_p|_{\{||x||\leq\epsilon\}}$ 
%is 
%Then by Lemma $\heartsuit\clubsuit\heartsuit\clubsuit$, 
\owari

\section{Construction of the cells}

\label{sec:constcell}

\subsection{Stable manifolds for real arrangements}

\label{subsec:stable}

The Lefschetz Theorem (\ref{thm:zale}) asserts that $\cpl(Q)$ has the 
homotopy type of a space obtained from $\cpl(Q)\cap F_\bbC$ by attaching 
some $\ell$-cells. The homotopy types of the attached cells 
are characterized by Theorem \ref{thm:char} 
under the Morse-Smale condition. 
In the remainder of this paper, we investigate the complexified 
real case, i.e., where each hyperplane $H\in\calA$ is 
defined by a linear equation with real coefficients. 
Let us briefly recall the set-up. 

Let $\calA=\{H_1, \ldots, H_n, H_\infty\}$ be an essential 
hyperplane arrangement in $\bbP_\bbC^\ell$, and 
$\alpha_i$ be the defining linear form of $H_i$ which is 
assumed to have real coefficients. 
Let $F=\{f=0\}$ be a generic hyperplane defined by a real linear 
form $f$. 
From Theorem \ref{thm:otvar}, there exist positive 
even integers 
$\lambda_0, \lambda_1, \ldots, \lambda_n, \lambda_\infty\in 2\bbZ_{>0}$ 
such that 
$\lambda_0=\lambda_1+\lambda_2+\cdots +\lambda_n+\lambda_\infty$ and 
$$
\varphi=\frac{f^{\lambda_0}}{\alpha_1^{\lambda_1}\cdots\alpha_n^{\lambda_n}
\alpha_\infty^{\lambda_\infty}}
$$
has only nondegenerate isolated critical points. 

The space 
$\bbP_\bbC^\ell-H_\infty$ is isomorphic to the affine space 
$\bbC^\ell$. We also denote by $\calA$ the induced 
affine arrangement $\{H_1, \ldots, H_n\}$ in $\bbC^\ell$ 
and by $\alpha_i$ the defining equation ($\deg =1$, with 
real coefficients) of $H_i$.  
Let $\sfch(\calA)$ be the set of all chambers of 
$\calA\cap\bbR^\ell$ and 
$\sfch_\ell^F(\calA)$ be the set of all chambers which 
do not meet $F_\bbR$. 
Denote by $\sfM(\calA)$ the complexified complement 
$\bbC^\ell-\bigcup_{i=1}^nH_i$. 

Let $C\in\sfch_\ell^F(\calA)$. Then 
$\varphi|_C$ is a positive real valued function and it 
has poles along the boundary $\partial \bar{C}$. 
Hence, for each $C$, $\varphi|_C$ has 
at least one critical point $p_C\in\rmint (C)$ in the relative 
interior of $C$. 
Then it follows from 
the Cauchy-Riemann equation that 
$p_C\in\sfM(\calA)$ is indeed a critical point of the function 
$\varphi:\sfM(\calA)\rightarrow\bbC$. 
Thus we obtain $|\sfch_\ell^F(\calA)|$ many critical points. 
From the assumption, $\varphi$ has only nondegenerate 
isolated critical points, the number of which is 
the Euler characteristic $|\chi(\sfM(\calA)-F_\bbC)|=b_\ell(\sfM(\calA))$
(see Proposition \ref{prop:flag}). 

\begin{Prop}
For each chamber $C\in\sfch_\ell^F(\calA)$ which does not 
meet $F_\bbR$, there exists only one critical point $p_C\in C$ of 
$\varphi$ in $C$. Conversely, any critical point is obtained 
in this way. 
\end{Prop}
In other words, the set of critical points $\rmCrit(\varphi)$ 
is parametrized by $\sfch_\ell^F(\calA)$. 
Moreover, since 
$|\varphi(z_1, \ldots, z_\ell)|=|\varphi(\bar{z}_1, \ldots, \bar{z}_\ell)|$, 
the gradient vector field $-\grad |\varphi|$ is invariant 
under complex conjugation. Thus we have the following: 

\begin{Thm}
\label{thm:stable}
The stable manifold of the critical point $p_C$ corresponding to 
a chamber $C\in\sfch_\ell^F(\calA)$ is 
$W_{p_C}^s=C\subset\sfM(\calA)$. 
\end{Thm}
In particular, the closure of a chamber $C\in\sfch_\ell^F(\calA)$ 
contains only one critical point $p_C$. Thus we have the following 
result. 

\begin{Cor}
The function $\varphi$ satisfies the Morse-Smale condition. 
\end{Cor}

Let $C\in\sfch_\ell^\calF(\calA)$ and 
$p_C\in C$ the corresponding critical point. We denote the 
attaching map of the unstable manifold $W_{p_C}^u$ by 
$\sigma_C:(\rmD^\ell, \partial\rmD^\ell)
\rightarrow(\sfM(\calA), \sfM(\calA)\cap F_\bbC)$. 
Since 
$H_\ell(\sfM(\calA), \bbC)
\stackrel{\sim}{\longrightarrow}
H_\ell(\sfM(\calA), \sfM(\calA)\cap F_\bbC; \bbC)$ is an isomorphism, 
$[\sigma_C]$ can be considered as an element of 
$H_\ell(\sfM(\calA), \bbC)$. Moreover 
$\{[\sigma_C]\}_{C\in\sfch_\ell(\calA)}$ forms a 
basis of $H_\ell(\sfM(\calA), \bbC)$. 
By Poincar\'e duality we have the following result. 

\begin{Cor}
\label{cor:basis}
$\{[C]\}_{C\in\sfch_\ell(\calA)}
\subset H_\ell^{lf}(\sfM(\calA), \bbC)$ 
forms a basis, and under suitable orientations, it is the 
dual basis of $\{[\sigma_C]\}_{C\in\sfch_\ell(\calA)}$. 
%\cong 
%H_\ell(\sfM(\calA), \bbC)$ and 
%$[C]\in H_\ell^{lf}(\sfM(\calA), \bbC)$ form a dual 
%basis. 
\end{Cor}

%\noindent
%{\bf Proof. } 
%\owari

Combining Theorem \ref{thm:stable} with Theorem \ref{thm:diffeo}, we 
easily prove the following result which is well known for 
$\ell=1$ (Example \ref{ex:1dim}). 

\begin{Cor}
$\sfM(\calA)\setminus \bigcup_{C\in\sfch_\ell(\calA)}C$ is 
diffeomorphic to $(F_\bbC-\calA)\times \rmD^2$. 
\end{Cor}

\begin{Example}
\label{ex:1dim}
\normalfont
Assume $\ell=1$. Let $\{a_1, a_2,\ldots, a_n\}\subset\bbR$ be an 
arrangement in $\bbR$ (we assume $a_1<\cdots<a_n$). Take a 
generic hyperplane (in this case, just a point) $F\in\bbR$ such that 
$a_i<F<a_{i+1}$. Then the set of chambers which do not meet 
$F$ is 
$$
\sfch_\ell(\calA)=\{(-\infty, a_1), \ldots, (a_{i-1}, a_i), 
(a_{i+1}, a_{i+2}), \ldots, (a_n, \infty)\}. 
$$
Hence 
$$
\sfM(\calA)\setminus \bigcup_{C\in\sfch_\ell(\calA)}C=
\bbC-\left([-\infty, a_i)\cup [a_{i+1}, \infty)\right), 
$$
which is diffeomorphic to $\rmD^2$. 
\end{Example}

\subsection{Construction of the Cells}
\label{subsec:const}

Our next task is to construct the cells 
in $\sfM(\calA)$ explicitly. More precisely, for a chamber 
$C\in\sfch_\ell^F(\calA)$ with $C\cap F_\bbR=\emptyset$ and 
fixed $p\in C$, we construct 
a continuous map 
$\sigma_C:(\rmD^\ell, \partial\rmD^\ell)\rightarrow 
(\sfM(\calA), \sfM(\calA)\cap F_\bbC)$ 
which is differentiable in a neighborhood of $0\in\rmD^\ell$, such that 
(recall the conditions in \S\ref{subsec:char}) 
\begin{itemize}
\item[{\tt (i)}] $\sigma_C(0)=p$, $\sigma_C(\rmD^\ell)\cap C=\{p\}$ and 
$\sigma_C(\rmD^\ell)$ intersects $C$ transversally. 
\item[{\tt (ii)}] $\sigma_C(\partial \rmD^\ell)\subset \sfM(\calA)\cap F_\bbC$. 
\item[{\tt (iii)}] If $C'\in\sfch_\ell^F(\calA)$ is another chamber, then 
$\sigma_C(\rmD^\ell)\cap C'=\emptyset$. 
\end{itemize}
Let us choose coordinates $(x_1, x_2, \ldots, x_\ell)$ such that 
$F$ is defined by $\{x_\ell=0\}$ and $p$ is 
$(0, 0, \ldots, 0,1)$. 
Recall that $L^{\ell-1}(\calA)$ is the set of all 
one-dimensional intersections of $\calA$. 
We can find a wide cylinder of height $1$ 
which ties up affine lines $L^{\ell-1}(\calA)$. 
More precisely, since $F$ is generic, each line 
$X\in L^{\ell-1}(\calA)$ intersects $F_\bbR$ 
transversely. Hence 
$$
R=2\sup\left\{\left.
\sqrt{x_1^2+\cdots +x_{\ell-1}^2}\right|\ 
(x_1, \ldots, x_{\ell-1}, x_\ell)
\in X\in L^{\ell-1}(\calA), 
0\leq x_\ell\leq 1\right\}
$$
is finite. Consider the cylinder 
$$
\rmCyl(R)=\{(x_1, \ldots, x_\ell)\in\bbR^\ell| 
x_1^2+\cdots +x_{\ell-1}^2=R^2, 0\leq x_\ell\leq 1\} 
$$
of radius $R$ and height $1$. 
The cylinder $\rmCyl(R)$ is diffeomorphic to 
$S^{\ell-2}(1)\times [0,1]$ under the map 
$$
\begin{array}{ccc}
S^{\ell-2}(1)\times [0,1]&\longrightarrow &\rmCyl(R)\\
&&\\
(x', t)&\longmapsto&(Rx', t), 
\end{array}
$$
where 
$x'=(x_1, \ldots, x_{\ell-1})\in S^{\ell-2}(1)=
\{(x_1, \ldots, x_{\ell-1});\ x_1^2+\cdots +x_{\ell-1}^2=1\}$. 
The boundary $\partial\rmCyl(R)$ of 
$\rmCyl(R)$ is the disjoint union of two spheres 
$S_0$ and $S_1$, where $S_t$ is a horizontal $(\ell-1)$-dimensional 
sphere of radius $R$ 
\begin{equation}
\label{eq:sphere}
\begin{array}{l}
S_t=\{(x_1, \ldots, x_\ell)|x_1^2+\cdots +x_{\ell-1}^2=R^2, x_\ell=t\}\\ 
D_t=\{(x_1, \ldots, x_\ell)|x_1^2+\cdots +x_{\ell-1}^2\leq R^2, x_\ell=t\}. 
\end{array}
\end{equation}

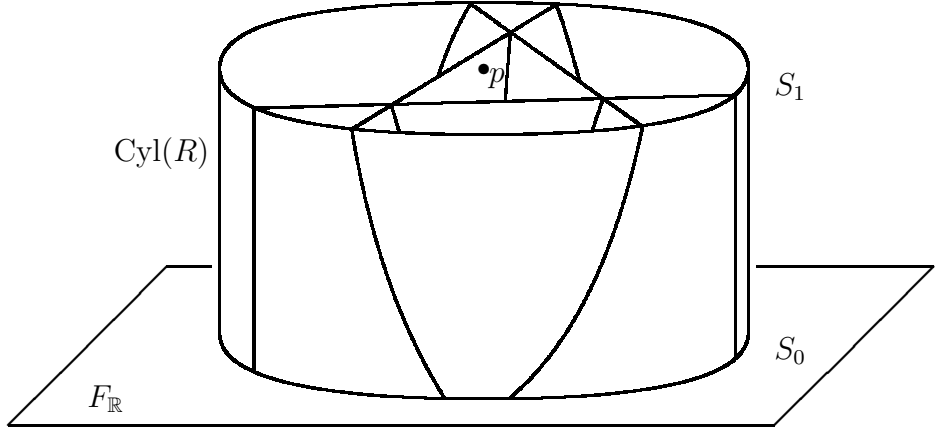
\begin{figure}[htbp]
\begin{picture}(100,150)(0,0)
\thicklines

\put(200,125){\circle*{4}}
\put(202,120){$p$}
\qbezier(100,125)(100,150)(200,150)
\qbezier(100,125)(100,100)(200,100)
\qbezier(300,125)(300,150)(200,150)
\qbezier(300,125)(300,100)(200,100)

\qbezier(100,25)(100,0)(200,0)
\qbezier(300,25)(300,0)(200,0)
\multiput(100,25)(200,0){2}{\line(0,1){100}}

\qbezier(113,110)(113,110)(295,115)
\put(113,110){\line(0,-1){100}}
\put(295,115){\line(0,-1){100}}

\qbezier(150,102)(150,102)(227,149)
\qbezier(150,102)(160,40)(185,0.5)

\qbezier(260,103)(260,103)(195,149.5)
\qbezier(260,103)(245,30)(210,0.5)

\put(20,-10){\line(1,1){60}}
\put(80,50){\line(1,0){17}}
\put(20,-10){\line(1,0){290}}
\put(310,-10){\line(1,1){60}}
\put(370,50){\line(-1,0){67}}

\put(310,15){$S_0$}
\put(310,115){$S_1$}
\put(50,-3){$F_\bbR$}
\put(60,90){$\rmCyl(R)$}

\qbezier(210,138)(210,138)(208,112.5)
\qbezier(165,111)(165,111)(168,101)
\qbezier(245,113.5)(245,113.5)(241,101.5)
\qbezier(195,149.5)(188,138.25)(183,123)
\qbezier(227,149)(228,150)(236,121)

\end{picture}
     \caption{The cylinder}\label{fig:cyl}
\end{figure}

$S_t\cap\calA$ determines a hypersphere arrangement 
on $S_t$ for $0\leq t\leq 1$. Note that the combinatorial 
type of this arrangement is independent of $t$. 

Since every $H\in\calA$ contains at least 
one affine line $X\in L^{\ell-1}(\calA)$ and $X$ intersects 
both $D_0$ and $D_1$, there exists a 
non-horizontal tangent vector $V_q\in\rmT_q(\rmCyl(R)\cap H)$ 
at each $q\in\rmCyl(R)\cap H$. 
Using a partition of unity, we have a vector field 
$\tilde{V}$ on $\rmCyl(R)$ of the form 
$$
\tilde{V}=\frac{\partial}{\partial x_\ell}+
\sum_{i=1}^{\ell-1}f_i\frac{\partial}{\partial x_i}
$$
which is tangent to each hypercylinder $\rmCyl(R)\cap H$. 
Now consider the one-parameter flow generated by $\tilde{V}$. 
The flow determines a diffeomorphism 
$\eta_t:S^{\ell-2}(1)\rightarrow S^{\ell-2}(1)$ such that 
$\eta_0(x')=x'$ and 
$$
\frac{d}{dt}(R\cdot \eta_t(x'), t)=\tilde{V}_{(R\cdot\eta_t(x'), t)}, 
$$
for $0\leq t\leq 1$ and $x'\in S^{\ell-2}(1)$. 
This determines a diffeomorphism 
$\eta_t:(S_0, S_0\cap\calA)\stackrel{\cong}{\longrightarrow}
(S_t, S_t\cap\calA)$. 

Let $\iota:S_1\rightarrow S_1$ be the involution 
$\iota :(x_1, \ldots, x_{\ell-1}, 1)\mapsto
(-x_1, \ldots, -x_{\ell-1}, 1)$ and define 
$I=\eta_1^{-1}\circ\iota\circ\eta_1:S_0\rightarrow S_0$. 
The next lemma follows immediately from the construction. 

\begin{Lemma}
\label{lem:equator}
Suppose $q\in S_0\cap H$ with $H\in\calA$, then the vectors 
$\overrightarrow{q I(q)}, \overrightarrow{qp}\in\rmT_qV_\bbR$ 
are on the same side with respect to the hyperplane $H\subset \rmT_qV_\bbR$. 
The same holds for $\overrightarrow{I(q)q}, \overrightarrow{I(q)p}\in
\rmT_{I(q)}V_\bbR$ 
when $I(q)\in S_0\cap H$. 
\end{Lemma}

Before defining 
$\sigma_C:(\rmD^\ell, \partial\rmD^\ell)
\rightarrow (\sfM(\calA), \sfM(\calA)\cap F_\bbR)$, 
we decompose the disk $\rmD^\ell$ into four pieces. Denote the latitude 
of $v\in\rmD^\ell$ by $\theta$, i.e., 
$v=(x'\cos\theta , ||x'||\sin\theta)$, where 
$x'=(x_1, \ldots, x_{\ell-1})$ with 
$x_1^2+\cdots +x_{\ell-1}^2\leq 1$. Fix 
$0<\theta_0<\pi/2$ so that $\tan\theta_0 =\frac{1}{R}$. 

\begin{itemize}
\item[(1)] (The core): $A_1=\left\{v\in\rmD^\ell\left| ||v||\leq \frac{1}{2}
\right.\right\}$. 
\item[(2)] (The northern hemisphere): 
$A_2=\left\{v\in\rmD^\ell\left| \frac{1}{2}\leq ||v||\leq 1, 
\theta\geq \theta_0\right.\right\}$. 
\item[(3)] (The southern hemisphere): 
$A_3=\left\{v\in\rmD^\ell\left| \frac{1}{2}\leq ||v||\leq 1, 
\theta\leq -\theta_0\right.\right\}$. 
\item[(4)] (The low latitudes): 
$A_4=\left\{v\in\rmD^\ell\left| \frac{1}{2}\leq ||v||\leq 1, 
-\theta_0\leq \theta\leq\theta_0\right.\right\}$. 
\end{itemize}

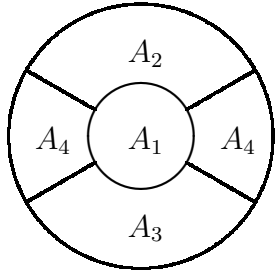
\begin{figure}[htbp]
\begin{picture}(100,100)(0,0)
\thicklines

\put(200,50){\circle{40}}

\qbezier(250,50)(250,63.39)(243.3,75)
\qbezier(243.3,75)(236.6,86.6)(225,93.3)
\qbezier(225,93.3)(213.39,100)(200,100)

\qbezier(175,93.3)(186.61,100)(200,100)
\qbezier(156.7,75)(163.4,86.6)(175,93.3)
\qbezier(150,50)(150,63.39)(156.7,75)

\qbezier(250,50)(250,36.61)(243.3,25)
\qbezier(243.3,25)(236.6,13.4)(225,6.7)
\qbezier(225,6.7)(213.39,0)(200,0)

\qbezier(175,6.7)(186.61,0)(200,0)
\qbezier(156.7,25)(163.4,13.4)(175,6.7)
\qbezier(150,50)(150,36.61)(156.7,25)

\qbezier(243.3,75)(243.3,75)(217.32,60)
\qbezier(243.3,25)(243.3,25)(217.32,40)
\qbezier(156.6,75)(156.6,75)(182.68,60)
\qbezier(156.6,25)(156.6,25)(182.68,40)

\put(195,45){$A_1$}
\put(195,78){$A_2$}
\put(195,12){$A_3$}
\put(230,45){$A_4$}
\put(160,45){$A_4$}

\end{picture}
     \caption{Decomposition of $\rmD^\ell$}\label{fig:dec}
\end{figure}

Given $v=(x'\cos\theta, ||x'||\sin\theta)\in\rmD^\ell$ 
with $1/2\leq ||x'||\leq 1$ 
and $\theta\neq 0$, let us define $\xi(v)\in V_\bbR$ by 
\begin{equation}
\label{eq:xi}
\xi(v):=
\left(
(2||x'||-1)\frac{\cos\theta}{\sin\theta}\cdot
\frac{-x'}{||x'||},\ 
2-2||x'||
\right).
\end{equation}
We also give an alternative description of $\xi(v)$. 
Straightforward computation shows that the line 
$p+t\cdot v$ ($t\in\bbR$) intersects the hyperplane $F_\bbR$ 
at $q=(-x'/(||x'||\tan\theta), 0)$. The point 
$\xi(v)$ above divide the segment $pq$ internally by the ratio 
$p\xi(v):\xi(v)q=(||x'||-\frac{1}{2}):(1-||x'||)$. 

We will define $\sigma_i:A_i\rightarrow \sfM(\calA)$, $i=1,\ldots,4$, 
separately. 
We use the notation in \S\ref{sec:cpxf} (\ref{eq:tg}) to express points 
in the complexified space $V_\bbC$. 
\begin{itemize}
\item[(1)]
$\sigma_1(v)=(p, v)_\bbC=p+\ii v\in V_\bbC$. 
%, in other words, $\sigma_1(v)$ 
%is represented by a tangent vector $v\in\rmT_pV_\bbR$ for $v\in A_1$. 
\item[(2)]
$\sigma_2(v)=(\xi(v), v)_\bbC$ 
for $v\in A_2$. 
The point $\sigma_2(v)$ is indeed contained in $\sfM(\calA)$ because every 
straight line passing through 
$p$ intersects each hyperplane $H\in\calA$ transversely. 
\item[(3)]
$\sigma_3(v)=(\xi(v), v)_\bbC$ 
for $v\in A_3$. 
\end{itemize}

The definition of $\sigma_4$ on the low latitudes $A_4$ is 
somewhat complicated. 
Let us define an annulus $T$ by 
$$
T=\{x'=(x_1, \ldots, x_{\ell-1});\ 1/2\leq ||x'||\leq 1\}. 
$$
Then the low latitudes $A_4$ can be expressed as 
$$
A_4=\{(x'\cos\theta, ||x'||\sin\theta)\in\rmD^\ell|\ 
x'\in T, -\theta_0\leq\theta\leq\theta_0\}. 
$$
We extend $\eta_t$ and $I:S^{\ell-2}(1)\rightarrow S^{\ell-2}(1)$ 
to $T$ by 
\begin{eqnarray*}
\eta_t(x')&:=&||x'||\cdot \eta_t\left(\frac{x'}{||x'||}\right), \\
I(x')&:=&\eta_1^{-1}(-\eta_1(x'))
\end{eqnarray*}
for $x'\in T$. Now define 
$$
\begin{array}{ccccc}
\gamma & : & T \times [0, 1]& \stackrel{\cong}{\longrightarrow}& A_4\\
&&&&\\
&&(x', t)&\longmapsto&
(-\eta_t^{-1}(-\eta_t(x'))\cos(2t-1)\theta_0, ||x'||\sin(2t-1)\theta_0). 
\end{array}
$$
Since $\eta_t:T\rightarrow T$ is a 
diffeomorphism, so is $\gamma$. 
Define 
$\sigma_4(\gamma(x',t))\in V_\bbC\cong\rmT V_\bbR$ by 
\begin{equation}
\label{eq:sigma4}
\sigma_4(\gamma(x',t))=
\left(
(1-t)\xi(\gamma(x',0))+t\xi(\gamma(x',1)), 
x'-I(x')
\right)_\bbC. 
\end{equation}

\begin{Lemma}
\label{lem:sigma4}
We have $\sigma_4(\gamma(x',t))\in\sfM(\calA)$. 
When $||x'||=1$, $\sigma_4(\gamma(x',t))$ is contained in 
$F_\bbC$, but is not contained in $F_\bbR$. 
\end{Lemma}

\noindent
{\bf Proof. } 
The second part is obvious. 
Indeed, since 
$x'-I(x')$ is a nonzero horizontal vector, it is 
contained in $\rmT F_\bbR$ when $||x'||=1$. 

Next we prove $\sigma_4(\gamma(x',t))\in\sfM(\calA)$ for 
$(x', t)\in T\times [0,1]$. 
By definition, we have 
\begin{eqnarray*}
&&\xi(\gamma(x',0))=
\left(
(2||x'||-1)
\frac{x'}{||x'||\tan\theta_0}, 
2-2||x'||
\right)\in V_\bbR, 
\\
&&\xi(\gamma(x',1))=
\left(
(2||x'||-1)
\frac{\eta_1^{-1}(-\eta_1(x'))}{||x'||\tan\theta_0}, 
2-2||x'||
\right)\in V_\bbR, 
\end{eqnarray*}
for $x'\in T$. 
Hence the tangent vector 
$x'-I(x')\in
\rmT_qV_\bbR$ with 
$q=(1-t)\xi(\gamma(x',0))+t\xi(\gamma(x',1))$, is parallel 
to $\xi(\gamma(x',0))-\xi(\gamma(x',1))$. 
In order to prove this lemma, it suffices to prove that 
the line segment connecting $\xi(\gamma(x',0))$ and $\xi(\gamma(x',1))$ 
is not contained in any hyperplane $H\in\calA$. 
So it suffices to prove the next lemma. 

%\owari

\begin{Lemma}
\label{lem:sameside}
If $\xi(\gamma(x',0))$ (resp. $\xi(\gamma(x',1))$) 
is contained in a hyperplane $H\in\calA$, 
then $\xi(\gamma(x',1))$  (resp. $\xi(\gamma(x',0))$) and 
$p$ lie on the same side with respect to $H$. 
\end{Lemma}

\noindent
{\bf Proof. } 
Suppose $\xi(\gamma(x',0))$ is contained in a 
hyperplane $H\in\calA$. 
Choose a defining equation $\alpha_H$ of $H$ 
such that $\alpha_H(p)>0$. We prove 
\begin{equation}
\label{eq:pos}
\alpha_H(\xi(\gamma(x',1)))>0. 
\end{equation}
Let us put 
$q=(x'/(||x'||\tan\theta_0), 0)=(Rx'/||x'||, 0)\in F_\bbR$. 
Since $\xi(\gamma(x',0))$ divides the segment $pq$ internally, 
we have $\alpha_H(q)<0$. 
By the definition of $\eta_t$, $(R\eta_t(x'/||x'||), t)\in V_\bbR$ 
($0\leq t\leq 1$) is a flow which is tangent to $H\cap\rmCyl_R$. 
Hence $\alpha_H(R\eta_t(x'/||x'||), t)<0$ for all $0\leq t\leq 1$. 
In particular, we have $\alpha_H(R\eta_1(x'/||x'||), 1)<0$. 
Since $p$ is the midpoint of $(R\eta_1(x'/||x'||), 1)$ and 
$(-R\eta_1(x'/||x'||), 1)$, we have $\alpha_H(-R\eta_1(x'/||x'||), 1)>0$. 
Similarly, using the flow $\eta_t$, we verify 
%Again using the flow $\eta_t$ to get 
$\alpha_H(R\eta_1^{-1}(-\eta_1(x'/||x'||), 0))>0$. 
$\xi(\gamma(x',1))$ divides the segment connecting 
$p$ and $(R\eta_1^{-1}(-\eta_1(x'/||x'||)), 0)$ internally, thus we have 
(\ref{eq:pos}). 

Similarly, if $\alpha_H(\xi(\gamma(x', 1)))=0$, then we have 
$\alpha_H(\sigma(\gamma(x', 0)))>0$. 
\owari

Now we are ready to construct the cell 
$\sigma_C:(\rmD^\ell, \partial\rmD^\ell)\rightarrow 
(\sfM(\calA), F_\bbC\cap\sfM(\calA))$. 
By definition, we have $\sigma_1|_{A_1\cap A_2}=\sigma_2|_{A_1\cap A_2}$ and 
$\sigma_1|_{A_1\cap A_3}=\sigma_3|_{A_1\cap A_3}$. Hence we have a continuous 
map 
$$
\sigma_{123}:A_1\cup A_2\cup A_3\longrightarrow \sfM(\calA). 
$$ 
Unfortunately, $\sigma_{123}$ and $\sigma_4$ do not coincide on their 
boundaries. However, we can paste the pieces together. 
Indeed, given a point 
$v=(x'\cos\theta_0, ||x'||\sin\theta_0)\in A_2\cap A_4$, 
both $\sigma_2(v)$ and $\sigma_4(v)$ can be considered 
as elements in $\rmT_{\xi(v)}V_\bbR$. 
Under the natural identification $\rmT_{\xi(v)}V_\bbR\cong V_\bbR$, 
we have 
\begin{eqnarray*}
&&\sigma_2(v)=v\\ 
&&\sigma_4(v)=x'-I(x'). 
\end{eqnarray*}
Recall that $v$ and $x'-I(x')$ are 
positive multiples 
of $p-\xi(v)=\overrightarrow{\xi(v)p}$ and 
$\xi(\gamma(x',0))-\xi(\gamma(x',1))$ respectively. 
Even if $\xi(v)$ is contained in some $H\in\calA$, 
$\sigma_2(v)$ and $\sigma_4(v)\in\rmT_{\xi(v)}V_\bbR\cong V_\bbR$ are 
on the same side with respect to a hyperplane $H\subset V_\bbR$ 
from Lemma \ref{lem:sameside}. We can continuously connect them by 
$$
\sigma_2(v)\cos\rho+\sigma_4(v)\sin\rho,\ 0\leq\rho\leq\frac{\pi}{2}, 
$$
thus we have a continuous map $\sigma_{1234}:\rmD^\ell\rightarrow\sfM(\calA)$ 
which satisfies {\tt (i)} and {\tt (iii)}. 
To glue the boundary $\partial\rmD^\ell$ to $F_\bbC$, we apply 
the following lemma. 
Let us set 
$$
\rmCyl_{R, \varepsilon}:=\{(x_1, \ldots, x_\ell)| x_1^2+\cdots+x_{\ell-1}^2\leq R^2, 
-\varepsilon\leq x_\ell\leq \varepsilon\}
$$
and 
$$
\rmCyl_{R, 0}:=\{(x_1, \ldots, x_\ell)| x_1^2+\cdots+x_{\ell-1}^2\leq R^2, 
x_\ell=0\}. 
$$

\begin{Lemma}
\label{lem:prod}
For sufficiently small $\varepsilon>0$, 
$(\rmCyl_{R,\varepsilon}, \rmCyl_{R,\varepsilon}\cap\calA)$ 
is diffeomorphic to 
$(\rmCyl_{R,0}, \rmCyl_{R,0}\cap\calA)\times[-\varepsilon, \varepsilon]$. 
\end{Lemma}
Denote the composite map 
$\rmCyl_{R,\varepsilon}\rightarrow\rmCyl_{R,0}\times[-\varepsilon, \varepsilon]
\rightarrow\rmCyl_{R,0}$ by $\rmPr_1$. 
Then, for $v\in\partial\rmD^\ell$, 
\begin{equation}
\label{eq:paste}
\sigma(v)\cos\rho+\rmPr_1(\sigma(v))\sin\rho,\ 0\leq\rho\leq\frac{\pi}{2}
\end{equation}
connects $\sigma(v)$ to $\rmPr_1(\sigma(v))\in\rmT_{\sigma(v)} F_\bbR$. 
Thus we have a 
map $\sigma_C:(\rmD^\ell, \partial\rmD^\ell)
\rightarrow(\sfM(\calA), \sfM(\calA)\cap F_\bbC)$ 
which satisfies {\tt (i)}, {\tt (ii)} and {\tt (iii)}. 
This completes the construction of the cell.

\begin{figure}[htbp]
\begin{picture}(100,520)(0,0)
\thicklines

%%%%%%%%
%%%%%%%%   e_1:A_1 --> M(A)
%%%%%%%%+390

\put(45,400){$A_1$}

\put(50,440){\circle{40}}
\put(50,440){\circle*{4}}
\put(50,440){\vector(1,2){8}}
\put(50,440){\vector(0,-1){12}}

\put(100,440){\vector(1,0){50}}
\put(123,448){$\sigma_1$}

\put(150,390){\line(1,0){200}}
\put(210,380){\line(1,1){110}}
\put(290,380){\line(-1,1){110}}
\put(310,465){$L_1$}
\put(185,465){$L_2$}
\put(340,395){$F_\bbR$}

\put(241,475){$p$}
\put(250,470){\circle*{4}}
\put(250,470){\circle{40}}
\put(250,470){\vector(1,2){8}}
\put(250,470){\vector(0,-1){12}}

%%%%%%%%
%%%%%%%%   e_2:A_2 --> M(A)
%%%%%%%%+260

\put(45,280){$A_2$}

\put(50,310){\circle*{4}}
\put(100,310){\vector(1,0){50}}
\put(123,318){$\sigma_2$}
\put(150,260){\line(1,0){200}}
\put(210,250){\line(1,1){110}}
\put(290,250){\line(-1,1){110}}
\put(230,340){$p$}
\put(250,340){\circle*{4}}

\qbezier(64.142,324.142)(57.654,330)(50,330)
\qbezier(35.858,324.142)(42.346,330)(50,330)
\qbezier(78.284,338.284)(65.308,350)(50,350)
\qbezier(21.716,338.284)(34.692,350)(50,350)
\qbezier(21.716,338.284)(21.716,338.284)(35.858,324.142)
\qbezier(64.142,324.142)(64.142,324.142)(78.284,338.284)

\put(50,310){\vector(1,4){5}}
\put(50,310){\vector(1,4){6.6}}
\put(50,310){\vector(1,4){8.3}}
\put(50,310){\vector(1,4){9.9}}

\put(250,340){\vector(1,4){4}}
\put(243.33,313.33){\vector(1,4){4.5}}
\put(243.33,313.33){\circle*{3}}
\put(236.66,286.66){\vector(1,4){5}}
\put(236.66,286.66){\circle*{3}}
\put(230,260){\vector(1,4){5.5}}
\put(230,260){\circle*{3}}

\put(250,340){\vector(1,1){12}}
\put(230,320){\vector(1,1){14}}
\put(230,320){\circle*{3}}
\put(210,300){\vector(1,1){16}}
\put(210,300){\circle*{3}}
\put(190,280){\vector(1,1){18}}
\put(190,280){\circle*{3}}
\put(170,260){\vector(1,1){19}}
\put(170,260){\circle*{3}}

\put(250,340){\vector(-1,1){12}}
\put(270,320){\vector(-1,1){14}}
\put(270,320){\circle*{3}}
\put(290,300){\vector(-1,1){16}}
\put(290,300){\circle*{3}}
\put(310,280){\vector(-1,1){18}}
\put(310,280){\circle*{3}}
\put(330,260){\vector(-1,1){19}}
\put(330,260){\circle*{3}}

\put(195,260){\vector(2,3){13}}
\put(195,260){\circle*{3}}

\put(220,260){\vector(1,2){10}}
\put(220,260){\circle*{3}}

\put(250,260){\vector(0,1){19}}
\put(250,260){\circle*{3}}

\put(280,260){\vector(-1,2){10}}
\put(280,260){\circle*{3}}

\put(305,260){\vector(-2,3){13}}
\put(305,260){\circle*{3}}

%%%%%%%%
%%%%%%%%   e_3:A_3 --> M(A)
%%%%%%%%+130

\put(45,120){$A_3$}

\put(50,180){\circle*{4}}
\put(100,180){\vector(1,0){50}}
\put(123,188){$\sigma_3$}
\put(150,130){\line(1,0){200}}
\put(210,120){\line(1,1){110}}
\put(290,120){\line(-1,1){110}}
\put(230,210){$p$}
\put(250,210){\circle*{4}}

\put(50,180){\vector(-1,-4){5}}

\qbezier(64.142,165.858)(57.654,160)(50,160)
\qbezier(35.858,165.858)(42.346,160)(50,160)
\qbezier(78.284,151.716)(65.308,140)(50,140)
\qbezier(21.716,151.716)(34.692,140)(50,140)
\qbezier(21.716,151.716)(21.716,151.716)(35.858,165.858)
\qbezier(64.142,165.858)(64.142,165.858)(78.284,151.716)

\put(250,210){\vector(-1,-1){12}}
\put(230,190){\vector(-1,-1){14}}
\put(230,190){\circle*{3}}
\put(210,170){\vector(-1,-1){16}}
\put(210,170){\circle*{3}}
\put(190,150){\vector(-1,-1){18}}
\put(190,150){\circle*{3}}

\put(250,210){\vector(-1,-4){4}}

\put(170,130){\vector(-1,-1){10}}
\put(170,130){\circle*{3}}

\put(250,210){\vector(1,-1){12}}
\put(270,190){\vector(1,-1){14}}
\put(270,190){\circle*{3}}
\put(290,170){\vector(1,-1){16}}
\put(290,170){\circle*{3}}
\put(310,150){\vector(1,-1){18}}
\put(310,150){\circle*{3}}

\put(330,130){\vector(1,-1){10}}
\put(330,130){\circle*{3}}

\put(195,130){\vector(-2,-3){8}}
\put(195,130){\circle*{3}}

\put(220,130){\vector(-1,-2){7}}
\put(220,130){\circle*{3}}

\put(250,130){\vector(0,-1){11}}
\put(250,130){\circle*{3}}

\put(280,130){\vector(1,-2){7}}
\put(280,130){\circle*{3}}

\put(305,130){\vector(2,-3){8}}
\put(305,130){\circle*{3}}

%%%%%%%%
%%%%%%%%   e_4:A_4 --> M(A)
%%%%%%%%

\put(45,10){$A_4$}

\put(50,50){\circle*{4}}
\put(100,50){\vector(1,0){50}}
\put(123,58){$\sigma_4$}
\put(150,0){\line(1,0){200}}
\put(210,-10){\line(1,1){110}}
\put(290,-10){\line(-1,1){110}}
\put(250,87){$p$}
\put(250,80){\circle*{4}}
\multiput(250,80)(4,-4){20}{\circle*{1.5}}
\multiput(250,80)(-4,-4){20}{\circle*{1.5}}

\put(50,50){\vector(-4,-1){18}}
\put(50,50){\vector(-3,2){25}}
\put(50,50){\vector(-3,-2){25}}
\put(50,50){\vector(-1,0){30}}

\qbezier(35.858,64.142)(30,57.654)(30,50)
\qbezier(35.858,35.858)(30,42.346)(30,50)
\qbezier(21.716,78.284)(10,65.308)(10,50)
\qbezier(21.716,21.716)(10,34.692)(10,50)
\qbezier(21.716,21.716)(21.716,21.716)(35.858,35.858)
\qbezier(35.858,64.142)(35.858,64.142)(21.716,78.284)

%\put(250,80){\vector(-1,-1){12}}
\put(230,60){\vector(-1,0){14}}
\put(230,60){\circle*{3}}
\put(210,40){\vector(-1,0){16}}
\put(210,40){\circle*{3}}
\put(190,20){\vector(-1,0){18}}
\put(190,20){\circle*{3}}

\put(250,80){\vector(-1,0){18}}

\put(170,2){\vector(-1,0){15}}
\put(170,1){\circle*{3}}

%\put(250,80){\vector(-1,1){12}}
\put(270,60){\vector(-1,0){14}}
\put(270,60){\circle*{3}}
\put(290,40){\vector(-1,0){15}}
\put(290,40){\circle*{3}}
\put(310,20){\vector(-1,0){16}}
\put(310,20){\circle*{3}}

\put(330,2){\vector(-1,0){17}}
\put(330,1){\circle*{3}}

\put(195,2){\vector(-1,0){18}}
\put(195,1){\circle*{3}}

\put(220,2){\vector(-1,0){16}}
\put(220,1){\circle*{3}}

\put(250,2){\vector(-1,0){15}}
\put(250,1){\circle*{3}}

\put(280,2){\vector(-1,0){16}}
\put(280,1){\circle*{3}}

\put(305,2){\vector(-1,0){18}}
\put(305,1){\circle*{3}}

\put(220,30){\vector(-1,0){18}}
\put(220,30){\circle*{3}}
\put(250,30){\vector(-1,0){20}}
\put(250,30){\circle*{3}}
\put(280,30){\vector(-1,0){18}}
\put(280,30){\circle*{3}}

\end{picture}
     \caption{$\sigma_i:A_i\rightarrow \sfM(\calA)$}\label{fig:ei}
\end{figure}

\begin{Example}
\normalfont
\label{ex:cell}
We illustrate the above construction for $\ell=2$. 
Let us consider an arrangement $\calA=\{L_1, L_2\}$ of 
two lines and a generic line $F$, 
\begin{eqnarray*}
L_1&:& y=x+\frac{1}{2}\\
L_2&:& y=-x+\frac{1}{2}\\
F&:& y=0.
\end{eqnarray*}
In this case, $\rmD^2$ is decomposed by the $A_i$ as in 
Figure \ref{fig:dec}. The map $\sigma_i:A_i\rightarrow\sfM(\calA)$ ($i=1,2,3,4$) 
is illustrated in Figure \ref{fig:ei}. 
\end{Example}

\section{Twisted minimal chain complexes}
\label{sec:twist}

The explicit construction in the previous section enables us 
to find an presentation of the cellular chain complex 
associated with the minimal CW decomposition 
with coefficients in a local system. We demonstrate 
this point in this section. 

\subsection{Flags and orientations}

In order to compute the boundary map of a cellular chain complex, 
we have to choose an orientation for each cell. 
First we recall some basic notation and terminology. 

Let $X$ be a differentiable manifold of $\dim_\bbR=n$ 
with boundary $\partial X$. Each orientation for $X$ 
determines an orientation for $\partial X$ as follows: 
Given $x\in\partial X$ choose a positively oriented 
basis $(v_1, v_2, \ldots, v_n)$ for $\rmT_xX$ in such a way 
that $v_2, \ldots, v_n\in\rmT_x(\partial X)$ and that 
$v_1$ is an outward vector. Then $(v_2, \ldots, v_n)$ 
determines an orientation on $\partial X$. 

Let $X$, $Y$ and $Z$ be oriented differentiable manifolds without boundary. 
Further assume $X$ is compact, $Z$ is a closed submanifold of $Y$, 
and $\dim X+\dim Z=\dim Y$. Let $f:X\rightarrow Y$ be 
differentiable map transversal 
to $Z$, i.e., 
$$
(df_x)(\rmT_xX)+\rmT_yZ=\rmT_yY
$$
holds at each point $x$ such that $y=f(x)\in Z$. 
Then $f^{-1}(Z)$ is a closed zero-dimensional submanifold 
of $X$, hence a finite set. 
Let $x\in f^{-1}(Z)$ and choose positively oriented bases 
$u=(u_1, \ldots, u_m)$ and $v=(v_1, \ldots, v_n)$ for 
$\rmT_xX$ and $\rmT_{f(x)}Z$, respectively. 
Under this assumption, we can define a local intersection number 
$I_x(f, Z)$ for each $x\in f^{-1}(Z)$ as follows: 
\begin{equation}
\label{eq:ori1}
I_x(f, Z)=
\left\{
\begin{array}{ll}
1&\mbox{ if } (f_*u, v) \mbox{ is positive for }\rmT_{f(x)}Y \\
-1&\mbox{ if } (f_*u, v) \mbox{ is negative for }\rmT_{f(x)}Y.  
\end{array}
\right. 
\end{equation}
And we also define $I(f, Z):=\sum_{x\in f^{-1}(Z)}I_x(f,Z)$. 

Let $V=\bbR^\ell$ be a real $\ell$-dimensional vector space 
and 
$$
\calF:\ \emptyset=\calF^{-1}\subset\calF^0\subset\calF^1\subset\cdots 
\subset\calF^\ell=V  
$$
be a complete flag of affine subspaces. 

\begin{Def}
\normalfont
An {\it oriented flag} is a flag $\calF$ in $V$ equipped with 
an orientation for each $\calF^i$, $i=1, \ldots, \ell$. 
\end{Def}
A given point $\calF^0\in V$ and a basis $v_1, \ldots, v_\ell$ 
of $V$ 
determine an oriented flag. Indeed, 
by defining 
$$
\calF^k:=\calF^0+\sum_{i=1}^{k}\bbR v_i, 
$$
$(v_1, \ldots, v_k)$ determines an orientation of 
$\calF^k$. Conversely, any oriented flag can be obtained 
in this way. 
Define positive and negative half subspaces, 
$\calF^k_+$ and $\calF^k_-$, by 
\begin{eqnarray*}
&&\calF^k_+=\calF^{k-1}+\bbR_{> 0}v_k\\
&&\calF^k_-=\calF^{k-1}+\bbR_{< 0}v_k, 
\end{eqnarray*}
respectively. 
Next we define signature of a chamber; a map 
$\sign :\sfch^k(\calA)\rightarrow \{\pm 1\}$. 

\begin{Def}
\normalfont
For $C\in\sfch^k(\calA)$, 
\begin{equation}
\sign(C)=
\left\{
\begin{array}{ll}
1&\mbox{ if }\calF^k\cap C\subset\calF^k_+\\
-1&\mbox{ if }\calF^k\cap C\subset\calF^k_-. 
\end{array}
\right.
\end{equation}
\end{Def}

\begin{Def}
\normalfont
Let $v_1, \ldots, v_\ell$ be a basis of $V_\bbR$. We fix an 
orientation of $V_\bbC$ by 
\begin{equation}
\label{eq:ori2}
(v_1, \ldots, v_\ell, \ii v_1, \ldots, \ii v_\ell). 
\end{equation}
Note that this orientation differs by a multiplication of 
$(-1)^{\ell(\ell-1)/2}$ from the canonical 
orientation defined by the complex structure 
$V_\bbC\cong\bbC^\ell$. 
\end{Def}

\subsection{Local systems and chambers}

Recall that 
$\sigma_C:(\rmD^\ell, \rmS^{\ell-1})\rightarrow 
(\sfM(\calA), \sfM(\calA)\cap F_\bbC)$ 
is the cell corresponding to the chamber $C\in\sfch_\ell^\calF(\calA)$. 
We can choose an orientation of $\sigma_C$ so that the intersection 
number satisfies 
\begin{equation}
\label{eq:ori3}
[C]\cdot [\sigma_C]=1. 
\end{equation}
Let $\Phi:\rmGal(\calA)\rightarrow \rmVect_\bbC$ be a representation 
of the Deligne groupoid and $\calL_\Phi$ be the associated 
local system. Since $C$ is a connected and 
simply connected subset of $\sfM(\calA)$, the space of 
flat sections $\calL_\Phi(C)$ is a finite dimensional vector space, 
and we have a natural isomorphism 
$$
\Phi(C)\cong \calL_\Phi(C). 
$$
From the fact that $\sfM(\calA)$ is homotopy equivalent to a space 
obtained from 
$\sfM(\calA)\cap F_\bbC$ by attaching 
$\ell$-cells $\{\sigma_C ;C\in\sfch_\ell(\calA)$\}, 
we have also a natural isomorphism 
\begin{equation}
\label{eq:chain}
H_\ell(\sfM(\calA), \sfM(\calA)\cap F_\bbC; \calL_\Phi), 
\cong
\bigoplus_{C\in\sfch_\ell(\calA)}\Phi(C)\otimes \bbC[\sigma_C], 
\end{equation}
where $\bbC[\sigma_C]$ is a one-dimensional vector space spanned by 
$[\sigma_C]$. 

\begin{Def}
Let $\calA$, $\calF^\bullet$ and $\Phi$ as above. 
Define 
\begin{equation}
\label{eq:tmcc}
\calC_k:=
\calC_k(\calA, \calF, \Phi)=
\bigoplus_{C\in\sfch_k^\calF(\calA)}\Phi(C)\otimes \bbC[\sigma_C]. 
\end{equation}
\end{Def}
From the general theory of cellular chain complexes, there exists a chain 
boundary map $\partial_\Phi:\calC_k\rightarrow\calC_{k-1}$ such that 
$$
H_k(\calC_\bullet, \partial_\Phi)\cong H_k(\sfM(\calA), \calL_\Phi). 
$$
We will give an formula for 
$\partial_\Phi:\calC_\bullet\rightarrow\calC_{\bullet-1}$. 

Let $\calL$ be a local system on $\sfM(\calA)$. 
Let $X$ be an oriented compact $\ell$-dimensional $C^\infty$-manifold 
with boundary $\partial X$, possibly $\partial X=\emptyset$, and 
$$
f :(X, \partial X)\longrightarrow (\sfM(\calA), \sfM(\calA)\cap F_\bbC) 
$$
be a smooth map. We denote by $f^*\calL$ the pull back of $\calL$ by $f$. 
Fix $x\in X$ and suppose $S\subset \sfM(\calA)$ is a 
connected and simply connected subset containing $f(x)$. Then 
there exists a natural isomorphism 
$$
f_{x, S}:(f^*\calL)(x)\stackrel{\sim}{\longrightarrow}
\calL(S).  
$$
Given a section $\alpha\in (f^*\calL)(X)$ we have a morphism 
from the constant sheaf $\bbC_X$ to $f^*\calL$ defined by 
$t\longmapsto t\cdot\alpha$, and it induces a homomorphism 
$$
\alpha\otimes\bullet :
H_\ell(X, \partial X; \bbC)\longrightarrow 
H_\ell(X, \partial X; f^*\calL). 
$$
Denote the image of the fundamental 
class $[X]\in H_\ell(X, \partial X; \bbC)$ by 
$\alpha\otimes [X]\in H_\ell(X, \partial X; f^*\calL)$. 
Hence we have 
$f_*(\alpha\otimes [X])\in H_\ell(\sfM(\calA),\sfM(\calA)\cap F_\bbC; \calL)$. 

Next we express $f_*(\alpha\otimes [X])$ by using the decomposition 
(\ref{eq:chain}). Recall (\ref{cor:basis}) that 
$\{[C]\}_{C\in\sfch_\ell(\calA)}\subset H_\ell^{lf}(\sfM(\calA), \bbC)$ 
is the dual basis to 
$\{[\sigma_C]\}_{C\in\sfch_\ell(\calA)}\subset 
H_\ell(\sfM(\calA), \bbC)$, i.e., for $C_1, C_2\in\sfch_\ell^\calF(\calA)$ 
$$
[C_1]\cdot[\sigma_{C_2}]=
\left\{
\begin{array}{ll}
1 & \mbox{if }C_1=C_2 \\
0 & \mbox{if }C_1\neq C_2.
\end{array}
\right.
$$
Thus we have the following lemma. 

\begin{Lemma}
\label{thm:formula}
Assume that 
$f^{-1}(C)$ is a finite set for each 
$C\in\sfch_\ell(\calA)$. 
Given a section $\alpha\in(f^*\calL)(X)$, 
$f_*(\alpha\otimes [X])\in H_\ell(\sfM(\calA), \sfM(\calA)\cap F_\bbC; \calL)$ 
is expressed as 
$$
f_*(\alpha\otimes [X])=
(-1)^\ell
\sum_{C\in\sfch_\ell^\calF(\calA)}
\sum_{x\in f^{-1}(C)}I_x(f, C)f_{x,C}(\alpha)\otimes[\sigma_C]. 
$$
\end{Lemma}

\subsection{The degree map}
\label{subsec:deg}

For the purpose of describing the boundary map of the chain 
complex (\ref{eq:tmcc}), we employ here an additional 
map, the degree map 
$$
\deg :\sfch_k^\calF(\calA)\times\sfch_{k-1}^\calF(\calA)
\longrightarrow \bbZ 
$$
defined below. For simplicity, 
we shall consider only the case $k=\ell$ and write 
$\calF^{\ell-1}=F$. 

Let $C\in\sfch_\ell^\calF(\calA)$ and $C'\in\sfch_{\ell-1}^\calF(\calA)$. 
Recall (\S\ref{subsec:const} (\ref{eq:sphere})) 
that $D_0$ is an $(\ell-1)$-dimensional large disk in $F_\bbR$ such that 
$C'\cap F_\bbR\in\sfch(\calA\cap F_\bbR)$ is a bounded chamber 
if and only if $C'\cap F_\bbR\subset D_0$. 

\begin{Def}
$$
\calP(C'):=\overline{C'}\cap D_0. 
$$
\end{Def}
In particular, $\calP(C')$ is equal to $\overline{C'}\cap F_\bbR$ 
if $C'\cap F_\bbR$ 
is a bounded chamber. The set $\calP(C')$ is, in any case, a 
convex closed subset of $F_\bbR$ with piecewise smooth boundary 
$\partial\calP(C')$. 
Next we consider vector fields on $\partial\calP(C')$ tangent 
to $F_\bbR$. 

\begin{Def}
\normalfont
Let $U\in\Gamma(\partial\calP(C'), \rmT F_\bbR|_{\partial\calP(C')})$ 
be a vector field on $\partial\calP(C')$ tangent to $F_\bbR$. Then 
$U$ is said to be {\it directing to} $C\in\sfch_\ell^\calF(\calA)$ if 
for any point $x\in\partial\calP(C')$ and hyperplane $H\in\calA$ with 
$x\in H$, $U(x)\notin\rmT_xH$ (in particular, $U(x)\neq 0$) and 
$U(x)$ and $C$ are in the same half-space with respect to $H$. 
Moreover if $x\in S_0=\partial D_0$, then $U(x)$ is 
assumed to be an inward vector. 
\end{Def}
Now we define the degree map as the degree of a certain Gauss map. 

\begin{Def}
Let $C\in\sfch_\ell^\calF(\calA)$ and  
$C'\in\sfch_{\ell-1}^\calF(\calA)$. Let $U$ be a vector field 
on $\partial\calP(C')$ directing to $C$. 
Then 
$$
\deg(C,C'):=\deg
\left(
\frac{U}{|U|}:\partial\calP(C')\longrightarrow S^{\ell-2}
\right). 
$$
\end{Def}
We need to prove the existence of a vector field 
$U\in\Gamma(\partial\calP(C'), \rmT F_\bbR|_{\partial\calP(C')})$ 
directing to $C$ and 
that $\deg(C,C')$ does not depend on the 
choice of $U$. 
From the genericity of $F$, there exists a tubular neighborhood 
$\calT\subset V_\bbR$ of $F_\bbR$ in $V_\bbR$ 
with a diffeomorphism (see also Lemma \ref{lem:prod})
$$
\tau:(\calT; \calA\cap\calT, \rmCyl_{R,\varepsilon}\cap\calT)
\stackrel{\sim}{\longrightarrow}
(F_\bbR; \calA\cap F_\bbR, S_0)\times (-1,1). 
$$
Fix a point $p\in C$. Then for $x\in\partial\calP(C')$, the vector 
$\overrightarrow{xp}\in\rmT_xV_\bbR$ is obviously in the 
same half-space 
as $C$ is for any $H\in\calA$ which contains $x$. 
The projection of this tangent vectors to $F_\bbR$ 
satisfies the condition, 
more precisely, 
\begin{equation}
\label{eq:vf}
U_\tau(x):=(\rmPr_1\circ\tau)_*(\overrightarrow{xp})
\in\rmT_xF_\bbR 
\end{equation}
determines a vector field on $\partial\calP(C')$ tangent to 
$F_\bbR$ directing to $C$, where $\rmPr_1$ is the first projection. 

Suppose $U$ and $U'$ are vector fields directing to $C$. 
Let $x\in\partial\calP(C')$. Consider the set $\calA_x$ 
of all hyperplanes in $\calA$ containing $x$. Then both 
$U(x)$ and $U'(x)\in\rmT_xF_\bbR$ are contained in the same 
chamber of $\calA_x$ which is also contains $C$. Hence $(1-t)U+tU'$ 
($0\leq t\leq 1$) is a 
continuous family of vector fields directing to $C$, and 
the maps $U/|U|$ and $U'/|U'|:\partial\calP(C')\rightarrow S^{\ell-2}$ 
are homotopic. Thus the degree $\deg(C,C')$ is well-defined. 

%\begin{Problem}
%Give a more combinatorial definition for the degree 
%map\footnote{Emanuele Delucchi has recently suggested an 
%alternative approach to define the degree map by 
%using matroid bundles, 
%it seems to be a possible future direction. }. 
%\end{Problem}

\subsection{The boundary map}

Recall that 
an arrangement $\calA$ with an oriented generic flag 
$\calF=\calF^\bullet$ and a representation 
$\Phi:\rmGal(\calA)\rightarrow\rmVect_\bbC$ of 
the Deligne groupoid determine a chain complex 
$(\calC_\bullet, \partial_\Phi)$ defined by 
$$
\calC_k:=
\calC_k(\calA, \calF, \Phi)=
\bigoplus_{C\in\sfch_k^\calF(\calA)}\Phi(C)\otimes \bbC[\sigma_C] 
$$
such that 
$H_k(\calC_\bullet, \partial_\Phi)\cong H_k(\sfM(\calA), \calL_\Phi)$. 
In this section we describe the boundary map $\partial_\Phi$ by using 
the degree map. 

\begin{Thm}
$\partial_\Phi:\calC_k\rightarrow\calC_{k-1}$ is expressed as follows: 
$$
\partial_\Phi(a\otimes[\sigma_C])=
-\sign (C)\times
\sum_{C'\in\sfch_{k-1}^\calF(\calA)}
\deg(C,C')\Delta_\Phi(C,C')(a)\otimes[\sigma_{C'}]. 
$$
(For $\Delta_\Phi(C, C')$, see Definition \ref{def:diff}.)
\end{Thm}

\noindent
{\bf Proof. }
We consider only the case where $k=\ell$. 
Recall that $\rmD^\ell=\{v\in\bbR^\ell; ||v||\leq 1\}$ and 
that the cell attaching map 
$\sigma_C:(\rmD^\ell, \partial\rmD^\ell)\rightarrow 
(\sfM(\calA), \sfM(\calA)\cap\calF^{\ell-1})$ was 
constructed in \S\ref{subsec:const}. 
The pull back $\sigma_C^*\calL_\Phi$ is canonically isomorphic 
to the trivial local system $\Phi(C)$. 
Since 
$\sigma_C^*\calL_\Phi|_{\partial\rmD^\ell}\cong(\partial\sigma_C)^*\calL_\Phi$, 
we have 
$$
a\otimes [\partial\rmD^\ell]
\in
H_{\ell-1}(\partial\rmD^\ell, (\partial\sigma_C)^*\calL_\Phi). 
$$
So we have to investigate the element 
$$
\begin{array}{rcc}
%&&H_{\ell-1}(\sfM(\calA)\cap\calF^{\ell-1}, \calL_\Phi)\\
%&&\downarrow\\
\sigma_{C *}(a\otimes [\partial\rmD^\ell])
&\in&H_{\ell-1}(\sfM(\calA)\cap\calF^{\ell-1}, 
\sfM(\calA)\cap\calF^{\ell-2}; \calL_\Phi)\\
&&||\wr\\
&&
\calC_{k-1}=
\bigoplus_{C'\in\sfch_{k-1}^\calF(\calA)}\Phi(C')\otimes \bbC[\sigma_{C'}]. 
\end{array}
$$
Here we recall some properties of the attaching map 
$\partial\sigma_C:\partial\rmD^\ell\rightarrow\sfM(\calA)\cap\calF^{\ell-1}$. 
First $\partial\rmD^\ell$ is divided into three parts 
$A_i':=A_i\cap\partial\rmD^\ell$, $i=2,3,4$ 
(see 
\S\ref{subsec:const} for the definitions of $A_1, \ldots, A_4$)
more precisely, 
\begin{itemize}
\item[(2)] %The Northern hemisphere: 
$A_2':=\{v\in\rmD^\ell; ||v||=1, \theta\geq\theta_0\}$, 
\item[(3)] %The Southern hemisphere: 
$A_3':=\{v\in\rmD^\ell; ||v||=1, \theta\leq -\theta_0\}$, 
\item[(4)] %The low latitudes: 
$A_4':=\{v\in\rmD^\ell; ||v||=1, 
-\theta_0\leq\theta\leq\theta_0\}$, 
\end{itemize}
where $\theta$ is the latitude of $v$, namely, 
$v=(x_1, \ldots, x_\ell)=(x'\cos\theta, ||x'||\sin\theta)$, 
and $\theta_0$ is a small fixed latitude. 
Write 
$(\partial\sigma_C)_i:=(\partial\sigma_C)|_{A_i'}
:A_i'\rightarrow F_\bbC$.

In view of Lemma \ref{thm:formula}, we have to count intersections 
of the map 
\begin{equation}
(\partial\sigma_C): A_2'\cup A_3'\cup A_4'\longrightarrow \calF_\bbC^{\ell-1} 
\end{equation}
with $C'\cap\calF^{\ell-1}_\bbR$ for 
$C'\in\sfch_{\ell-1}^\calF(\calA)$. 

It follows from Lemma \ref{lem:sigma4} that 
$(\partial\sigma_C)_4:
A_4'\rightarrow\sfM(\calA)\cap\calF_\bbR^{\ell-1}$ 
does not 
intersect $C'\cap\calF^{\ell-1}_\bbR$ for 
any chamber $C'\in\sfch_{\ell-1}^\calF(\calA)$. 

Suppose $v\in A'_2$. Recall that by the 
definition (\ref{eq:xi}) of $\xi$, 
$\xi(v)\in\calF_\bbR^{\ell-1}$ is the point such that the vector 
$\overrightarrow{\xi(v)p}$ is proportional to $v$. 
And it is obvious that the map 
$\xi_{A'_2}: A'_2\rightarrow D_0: v\mapsto \xi(v)$ is 
a diffeomorphism. The orientation on $A'_2$ is determined by 
(\ref{eq:ori3}). $\xi_{A'_2}$ is orientation preserving (resp. reversing) 
if $\sign(C)=1$ (resp. $\sign(C)=-1$). 
$(\partial\sigma_C)_2(v)$ can be expressed as 
$$
(\partial\sigma_C)_2(v)=U_\tau(\xi(v))\in\rmT_{\xi(v)}F_\bbR. 
$$
The vector field $U_\tau$ is not zero on $(\calA\cap D_0)\cup S_0$. 
Up to small perturbation, we may 
assume that the zero locus of 
$U_\tau$ consists of a finite number of points. 
Intersections of $(\partial\sigma_C)_2$ with 
$\calP(C')\subset\sfM(\calA)\cap\calF_\bbC^{\ell-1}$ can be 
thought of as the set of  singular points of the vector field $U_\tau$. 
Hence the sum of local intersection numbers is equal to the 
degree of the map from the boundary $\partial\calP(C')$ to 
the sphere $S^{\ell-2}$. Thus we have 
$$
I((\partial\sigma_C)_2, \calP(C'))=
(-1)^{\ell-1}\sign(C)\deg(C,C'). 
$$
Similarly, 
$$
(\partial\sigma_C)_3(v)=-U_\tau(\xi(v))\in\rmT_{\xi(v)}F_\bbR, 
$$
and we have 
$$
I((\partial\sigma_C)_3, \calP(C'))=
(-1)^{\ell}\sign(C)\deg(C,C'). 
$$
By the definition of $(\sigma_C)_2:A_2\rightarrow\sfM(\calA)$, 
$(\partial\sigma_C)_{2*}:\Phi(C)\rightarrow (\partial\sigma_C)_2^*\Phi(C')
\stackrel{\sim}{\longrightarrow}\Phi(C')$ is equal to 
$$
\Phi_{P^-(C, C')}:\Phi(C)\longrightarrow \Phi(C'). 
$$
Similarly, 
$(\partial\sigma_C)_{3*}:\Phi(C)\rightarrow\Phi(C')$ is equal to 
$$
\Phi_{P^+(C, C')}:\Phi(C)\longrightarrow \Phi(C'). 
$$
The proof is then completed by employing 
Lemma \ref{thm:formula}. 
\owari

\subsection{Examples}
\label{subsec:ex}

Let $\calA=\{H_1, \ldots, H_n\}$ be a hyperplane arrangement 
in $\bbR^\ell$. Fix a nonzero complex number $q_i\in\bbC^*$ for 
each $i=1, \ldots, n$. Then we can define a representation $\Phi$ 
of $\rmGal(\calA)$ as follows. First we put 
$$
\Phi(C)=\bbC[\sigma_C]
$$
for each $C\in\sfch(\calA)$. Given two chambers 
$C, C'\in\sfch(\calA)$, suppose that 
$\{H_{i_1}, \ldots, H_{i_k}\}$ is 
the set of all hyperplanes separating $C$ and $C'$. Then 
define 
$$
\Phi_{P^+(C, C')}: \Phi(C)\longrightarrow \Phi(C'),\ 
[\sigma_C]\longmapsto q_{i_1}q_{i_2}\cdots q_{i_k} [\sigma_{C'}]. 
$$
By the definition of $\rmGal(\calA)$ (\ref{def:gal} (4)) 
$\Phi$ determines a representation 
$\Phi:\rmGal(\calA)\rightarrow\rmVect_\bbC$, hence a 
rank one local system $\calL_\Phi$, such that 
the local monodromy around $H_i$ is $q_i^2$. 
Conversely any rank one local system can be obtained in this way. 

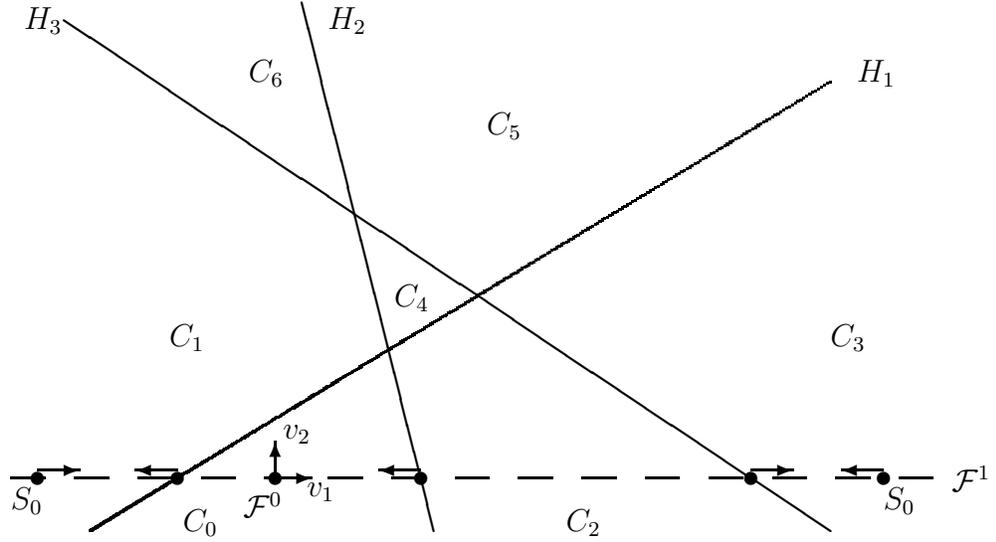
\begin{figure}[htbp]
\begin{picture}(100,200)(0,0)
\thicklines

\multiput(40,20)(24,0){15}{\line(1,0){13}}

\put(200,0){\line(-1,4){50}}
%\put(100,0){\line(4,3){230}}
\qbezier(70,0)(70,0)(350,170)
\put(350,0){\line(-3,2){290}}

\put(140,20){\circle*{5}}
\put(140,20){\vector(1,0){14}}
\put(140,20){\vector(0,1){14}}
\put(128,5){$\calF^0$}
\put(152,13){$v_1$}
\put(143,35){$v_2$}

\put(396,15){$\calF^1$}
\multiput(50,20)(320,0){2}{\circle*{5}}
\multiput(40,8)(330,0){2}{$S_0$}

\put(103,20){\circle*{5}}
\put(195,20){\circle*{5}}
\put(320,20){\circle*{5}}

\put(50,23){\vector(1,0){16}}
\put(103,23){\vector(-1,0){16}}
\put(195,23){\vector(-1,0){16}}
\put(320,23){\vector(1,0){16}}
\put(370,23){\vector(-1,0){16}}

\put(105,0){$C_0$}
\put(100,70){$C_1$}
\put(250,0){$C_2$}
\put(350,70){$C_3$}
\put(185,85){$C_4$}
\put(220,150){$C_5$}
\put(130,170){$C_6$}

\put(360,170){$H_1$}
\put(160,190){$H_2$}
\put(45,190){$H_3$}

\end{picture}
     \caption{Three lines with flags and a vector field directing 
to $C_6$}\label{fig:ex}
\end{figure}

Let us consider an arrangement of three lines $\calA$ (Fig. \ref{fig:ex}) 
with a generic flag 
$\calF^0\subset\calF^1\subset\calF^2=V_\bbR$ 
oriented by $(v_1, v_2)$ and 
\begin{eqnarray*}
\sfch_0^\calF(\calA)&=&\{C_0\}\\
\sfch_1^\calF(\calA)&=&\{C_1, C_2, C_3\}\\
\sfch_2^\calF(\calA)&=&\{C_4, C_5, C_6\}. 
\end{eqnarray*}
A vector field directing to $C_6$ is also drawn in Figure \ref{fig:ex}. 

Define the chain complex 
$\calC_k$ as in \S \ref{sec:twist}. 
The boundary map $\partial_\Phi:\calC_2\rightarrow\calC_1$ is, 
for example, 
\begin{eqnarray*}
\partial[\sigma_6]&=&
-\sum_{i=1}^3\deg(C_6, C_i)\Delta_\Phi(C_6, C_i)[\sigma_i]\\
&=&
-\left(
-(q_3-q_3^{-1})[\sigma_1]+(q_1q_2q_3-q_1^{-1}q_2^{-1}q_3^{-1})[\sigma_2]
-(q_1q_2-q_1^{-1}q_2^{-1})[\sigma_3]
\right). 
\end{eqnarray*}
Similarly, 
\begin{eqnarray*}
\partial[\sigma_5]&=&
(q_2q_3-q_2^{-1}q_3^{-1})[\sigma_1]+(q_1-q_1^{-1})[\sigma_3], \\
\partial[\sigma_4]&=&
(q_2-q_2^{-1})[\sigma_1]+(q_1-q_1^{-1})[\sigma_2]. 
\end{eqnarray*}
$\partial:\calC_1\rightarrow\calC_0$ is calculated as 
$$
\left(
\begin{array}{c}
\partial[\sigma_1]\\
\partial[\sigma_2]\\
\partial[\sigma_3]
\end{array}
\right)
=
\left(
\begin{array}{c}
-(q_1-q_1^{-1})\\
q_2-q_2^{-1}\\
q_2q_3-q_2^{-1}q_3^{-1}
\end{array}
\right)
[\sigma_0]. 
$$
Then direct computations show that the local system is resonant, i.e. 
$H_1(\calC_\bullet, \partial)\neq 0$, if and only if 
$q_1^2=q_2^2=q_3^2=1$. 

If we move the hyperplane $H_2$ so that the chamber $C_4$ collapses, 
we obtain another arrangement $\overline{\calA}$. In this case 
$\calC_2(\overline{\calA})$ is generated by $[\sigma_5]$ and $[\sigma_6]$. 
Hence the local system is resonant exactly when  
$\partial[\sigma_5]$ and $\partial[\sigma_6]$ are 
linearly dependent, or equivalently, $(q_1q_2q_3)^2=1$.

\section{Appendix}
\label{sec:fundam}

From attaching maps for $\ell=2$, 
we obtain a presentation for the fundamental group 
$\pi_1(\sfM(\calA))$. 

Let $\calA=\{H_1, \ldots, H_n\}$ be a line arrangement in $V=\bbR^2$. 
Let $\calF^0\subset\calF^1=F\subset V$ be an oriented generic flag. 
Note that $F_\bbR$ is an oriented line. 
We may assume that the chambers are ordered as 
$$
\begin{array}{ccl}
\sfch_0^\calF(\calA)&=&\{C_0\}\\
\sfch_1^\calF(\calA)&=&\{C_1, \ldots, C_n\}\\
\sfch_2^\calF(\calA)&=&\{C_{n+1}, \ldots, C_{n+b_2}\},
\end{array}
$$
and that the ordering $C_1, \ldots, C_n$ goes along with 
the orientation, that is, the intervals 
$C_1\cap F_\bbR, \ldots, C_n\cap F_\bbR$ are ordered from 
a negative place to a positive place with respect to 
the orientation for $F_\bbR$. 
The corresponding $1$-cells $\{\gamma_i=\sigma_{C_i}\}_{i=1, \ldots, n}$ 
are illustrated in Figure \ref{fig:loops}. 

\begin{figure}[htbp]
\begin{picture}(100,150)(0,0)
\thicklines

%\put(0,75){\vector(1,0){400}}

\put(50,75){\circle{4}}
\put(48,75){\line(-1,0){48}}
\put(52,75){\line(1,0){18}}
\multiput(75,75)(10,0){3}{\circle*{3}}
\put(123,75){\line(-1,0){18}}

\put(125,75){\circle{4}}

\put(160,75){\circle*{6}}
\put(155,55){$\calF^0$}

\put(200,75){\circle{4}}
\put(202,75){\line(1,0){46}}

\put(250,75){\circle{4}}
\put(252,75){\line(1,0){32}}
\put(348,75){\line(-1,0){23}}
\multiput(295,75)(10,0){3}{\circle*{3}}
\put(350,75){\circle{4}}
\put(352,75){\vector(1,0){50}}

\put(-2,60){$C_1$}
\put(60,60){$C_2$}
\put(205,60){$C_i$}
\put(355,60){$C_n$}

\put(400,60){$F_\bbR$}

\qbezier(160,75)(80,25)(50,25)
\put(50,25){\vector(-1,0){0}}
\qbezier(50,25)(20,25)(20,75)
\qbezier(160,75)(80,125)(50,125)
\put(55,125){\vector(1,0){0}}
\qbezier(50,125)(20,125)(20,75)
\put(20.1,85){\vector(0,1){0}}
\put(30,90){$\gamma_1$}

\qbezier(160,75)(200,50)(220,50)
\put(220,50){\vector(1,0){0}}
\qbezier(220,50)(235,50)(235,75)
\qbezier(160,75)(200,100)(220,100)
\put(216,100){\vector(-1,0){0}}
\qbezier(220,100)(235,100)(235,75)
\put(234.5,85){\vector(0,1){0}}
\put(240,90){$\gamma_i$}

\qbezier(160,75)(270,0)(340,0)
\put(340,0){\vector(1,0){0}}
\qbezier(340,0)(380,0)(380,75)
\qbezier(160,75)(270,150)(340,150)
\put(334,150){\vector(-1,0){0}}
\qbezier(340,150)(380,150)(380,75)
\put(379.8,85){\vector(0,1){0}}
\put(385,90){$\gamma_n$}

\end{picture}
     \caption{1-cells in $F_\bbR\otimes\bbC$}\label{fig:loops}
\end{figure}

Each chamber $C\in\sfch_2^\calF(\calA)$ is corresponding to a $2$-cell 
$\sigma_C$. Thus the boundary $\partial\sigma_C$, which is a 
word of generators $\{\gamma_i\}$, gives a relation in 
the fundamental group. The relation is 
$$
R(C):=
\ \gamma_1^{e_1}\gamma_2^{e_2}\cdots\gamma_n^{e_n}
\ \gamma_1^{-e_1}\gamma_2^{-e_2}\cdots\gamma_n^{-e_n}
=1, 
$$
where 
$$
e_i=\deg(C, C_i). 
$$

\begin{Thm}
The fundamental group $\pi_1(\sfM(\calA))$ is presented as: 
$$
\pi_1(\sfM(\calA))\cong
\left\langle
\gamma_1, \ldots, \gamma_n\ \left|\ R(C), C\in\sfch_2^\calF(\calA)
\right. \right\rangle. 
$$
\end{Thm}
We apply this theorem to the arrangement in Figure \ref{fig:ex} 
(\S\ref{subsec:ex}). Then 
$$
\begin{array}{ccl}
R(C_4)&=&\gamma_1^{-1}\gamma_2^{-1}\gamma_1\gamma_2\\
R(C_5)&=&\gamma_1^{-1}\gamma_3^{-1}\gamma_1\gamma_3\\
R(C_4)&=&\gamma_1^{-1}\gamma_2\gamma_3^{-1}
\gamma_1\gamma_2^{-1}\gamma_3. 
\end{array}
$$
Hence the fundamental group is isomorphic to the abelian group 
$$
\pi_1(\sfM(\calA))=\bbZ\gamma_1\oplus\bbZ\gamma_2\oplus\bbZ\gamma_3. 
$$

%kaiteru

%\section{}

\vspace{4mm}

\noindent
{\bf Acknowledgements. }
In the preparation of this paper, 
the author benefited from many people's advice. 
The author wishes to thank 
Manabu Akaho, 
David Bessis, 
Emanuele Delucchi, 
Yukihito Kawahara, 
Toshitake Kohno, 
Anatoly Libgober, 
Simona Settepanella, 
Hiroaki Terao, 
L\^e D\~ung Tr\'ang, 
Kyoji Saito and 
Sergey Yuzvinsky  
for their interests in this work, comments, helpful discussion and 
encouragements, 
Hiroyuki Yashima for discussing the Deligne groupoid. 
The author also thanks Eiko Fukunaga for her help on English. 
This research is supported by JSPS 
Fellowships for Young Scientists No. 16-0658.

\end{document}